\documentclass[11pt]{article}
\usepackage{graphicx,color,xcolor,enumitem,booktabs}
\usepackage[colorlinks=true, linkcolor=black, citecolor=black, urlcolor=black]{hyperref}
\usepackage{authblk}
\usepackage[normalem]{ulem}
\usepackage{geometry}
\usepackage[labelfont=bf]{caption}
\usepackage[english]{babel}
\usepackage{natbib}
\usepackage{bbm}

\usepackage[T1]{fontenc}       
\usepackage{lmodern}           
\usepackage[utf8]{inputenc}           
 \usepackage{amsmath, amssymb, amsfonts, amsthm} 
\usepackage{mathtools}        
\usepackage{mathrsfs}         
\usepackage{bbm}              

\newtheorem{theorem}{Theorem}[]
\newtheorem{lemma}[]{Lemma}
\newtheorem{proposition}[]{Proposition}

\theoremstyle{definition}

\newcommand{\one}{\mathbbm{1}}
\newcommand{\RR}{\mathbb{R}}
\newcommand{\DD}{\mathbb{D}}
\newcommand{\PP}{\mathbb{P}}
\newcommand{\EE}{\mathbb{E}}
\newcommand{\MM}{\mathbb{M}}
\newcommand{\GG}{\mathbb{G}}
\newcommand{\QQ}{\mathbb{Q}}

\begin{document}
\title{Deconvolution in unlinked linear models}
\date{\today}
\author[1]{Fadoua Balabdaoui}
\author[1,2]{Antonio Di Noia}
\author[3]{C\'ecile Durot}
\affil[1]{Seminar for Statistics, Department of Mathematics, ETH Zurich}
\affil[2]{Faculty of Economics, Euler Institute, Università della Svizzera italiana}
\affil[3]{MODAL'X, UMR CNRS 9023, Universit\'e Paris Nanterre}
\renewcommand\Affilfont{\itshape\small}
\maketitle

\begin{abstract}
Unlinked regression, in which covariates and responses are observed separately without known correspondence, has recently gained increasing attention.
Deconvolution, on the other hand, is a fundamental and challenging problem in nonparametric statistics with the aim of estimating the distribution of a latent random variable $Z$ based on observations contaminated by some additive noise. The complexity of this task is heavily influenced by the smoothness of the noise distribution and often leads to slow estimation rates.
In this paper, we combine the recent unlinked linear regression problem with the classical deconvolution framework. Specifically, we study nonparametric deconvolution under the assumption that $Z$ is a linear function of an observable multidimensional covariate. This structural constraint allows us to introduce a nonparametric estimator of the distribution of $Z$ which achieves the parametric rate of convergence in the Wasserstein distance of order 1, where the smoothness of the noise does not affect the rate.
Furthermore, we introduce nonparametric estimators for the unconditional density of $Z$ and the conditional density of $Z$ given an observed response. This allows us to study the problem of estimating the value of the latent linear predictor, whose link to the observed response is not accessible. Through several simulations, we illustrate the fast convergence rate of our deconvolution estimator and the performance of the proposed conditional estimators of the latent predictor in different simulation scenarios.
\end{abstract}


\section{Introduction}
\subsection{Background and contributions}
The unlinked or unmatched regression has recently regained attention in the statistics and machine learning literature. In this problem, we are given independent copies of covariates $X_1, \dots, X_n$ and responses $Y_1, \dots, Y_m$ with the same distribution as random variables $X$ and $Y$ respectively. Here, the sample sizes $n \ge 1$ and $m \ge 1$ are not necessarily equal. Assuming that $Y \overset{d}{=} f(X) + \epsilon$, the goal is to reconstruct the regression function $f$ under the condition that the noise $\epsilon$ has a known distribution. This problem arises in various fields, including sociology, economics, data privacy, computer vision, biology, and multi-tracking scenarios; see for example \cite{durot2024minimax}  and \cite{ slawski2024permuted} for a detailed account of the applications of unlinked regression. 

Theoretical developments in unlinked regression often leverage techniques from the deconvolution literature. For instance, the case where $f$ is monotone has been studied under various assumptions on noise distribution using deconvolution techniques \citep{carpentier2016learning, balabdaoui2021unlinked, rigollet2019uncoupled, meis}. In the unlinked linear regression setting, \cite{azadkia2022linear} proposed a deconvolution least squares estimator for parameter estimation, highlighting that identifiability issues may arise. Under some technical conditions on the noise and covariate distributions, consistency and asymptotic normality were established. The linear case has been previously addressed only in a closely related sub-problem, namely the \emph{permuted regression} or {\it shuffled regression}. In the latter setting, there exists a link between covariates and responses, however, this link is not known due to some unknown shuffling of the covariates (or equivalently the responses). The main focus in the existing literature on permuted regression is either exact or approximate recovery of the unknown permutation; 
see \cite{slawski2020two}, \cite{hsu2017linear}, \cite{pananjady2017linear}, \cite{zhang2021benefits}, \cite{slawski2021pseudo}, \cite{slawski2019linear}, \cite{tsakiris2020algebraic} and \cite{Unnikrishnan:2018gp}.
The gap between unlinked and permuted regression becomes particularly pronounced when the noise vanishes at a rate smaller than a threshold of order close to $n^{-1/2}$, as formalized in \cite{durot2024minimax}, where a phase transition in minimax rates under the Wasserstein distance was established: the minimax rates of estimation are the same in the two models when the noise level is larger than the threshold (in particular, under fixed noise level) and becomes smaller in the permuted regression than in unlinked regression when the noise level is of smaller order than $n^{-1/2}$. 

Deconvolution, on the other hand, is an old fundamental inverse problem in nonparametric statistics that arises in a wide range of applied fields, including econometrics, biometrics, medical statistics, image analysis, and signal processing. The main objective of deconvolution is to estimate an unknown density when only observations contaminated with additive noise are available. Note that it is often the case that the noise distribution is assumed to be known. Formally, given a latent random variable $Z$, we seek to estimate its distribution using observed data from $Y = Z + \epsilon$, where $\epsilon$ represents the additive noise, assumed to be independent of $Z$. The estimation task requires inverting the convolution induced by $\epsilon$ to recover the unknown distribution of $Z$.

Despite its seemingly straightforward formulation, the deconvolution problem is notably challenging. Its complexity, often measured in terms of the convergence rates in a chosen metric (e.g., $L^p$ or Wasserstein distance), is driven by the smoothness properties of the noise distribution. \cite{carroll1988optimal} established optimal pointwise estimation rates in the Gaussian case showing that they are very slow, specifically, a power of the logarithm of the sample size. The seminal papers \cite{fan1991optimal,fan1992deconvolution,fan1993adaptively} generalized this finding by showing that in the case where the characteristic function of the noise $\epsilon$ decays exponentially (as for the Gaussian noise), known as the \emph{supersmooth} case, the pointwise and global optimal rates in the $L^p$-norm are also a power of the logarithm of the sample size. In \cite{fan1991optimal,fan1993adaptively}, it is shown that for the \emph{ordinary smooth} case, i.e.\ where the characteristic function decays polynomially, the pointwise and global optimal rates in the $L^p$-norm are faster, specifically, a power of the sample size.

The statistical literature on deconvolution is very extensive. Here,  we shall provide only a few references thereof, and refer the reader to the book \cite{meister2009} for a nice overview. A variety of methodologies have been proposed in deconvolution problems, including kernel-based estimators which are based on the empirical characteristic function and Fourier inversion \citep{liu1989consistent, stefanski1990deconvolving, Fan91AsympNorm, fan1991optimal, es2005asymptotic}. Extensions to settings where the noise distribution is unknown, partially or fully, have been explored e.g. in \cite{butucea2005minimax,meister2006density, delaigle2014parametrically}.
Alternative estimation techniques include series estimators \citep{pensky1999adaptive, lounici2011global, carrasco2011spectral, carroll2004low, hall2005discrete}, nonparametric maximum likelihood methods \citep{groeneboom2003density, guan2021fast}, and penalization approaches \citep{comte2006penalized, comte2007finite}.

Our contribution in this paper bridges nonparametric deconvolution with the unlinked linear regression. Specifically, we focus on estimating the distribution of $Z \stackrel{d}{=} \beta^\top _0 X$ (for some non-unique $\beta_0$) from realizations of $Z + \epsilon$ and of the covariate $X$ whose distribution is unknown. This framework extends the classical deconvolution problem through incorporating the structure of the unlinked linear model. We develop an estimator based on the deconvolution least squares approach from \cite{azadkia2022linear} and the empirical measure of $X_1, \dots, X_n$. Using empirical process techniques, we show that this additional structure allows us to obtain the parametric rate of convergence in the Wasserstein distance for estimating the distribution of the latent variable $Z$, where the smoothness of the noise does not affect the rate.
We emphasize that the Wasserstein distance is a natural risk measure to consider for estimating a distribution, especially in cases where no regularity assumption is made. See e.g.\ \cite{caillerie2013deconvolution,dedecker2013minimax,dedecker2015improved}. In contrast, risk measures based for example on some $L^p$ distance between probability density functions can only be used for distributions that are assumed to have a probability density. Such assumptions are not required in this paper for the distribution of $Z$ to achieve the fast convergence rate in the sense of the Wasserstein risk. It is of crucial importance to emphasize that although the model incorporates a finite-dimensional parameter, it does not fall within the framework of semi-parametric models. In fact, our goal is not the estimate the finite-dimensional regression vector  $\beta_0$ but rather an infinite-dimensional object: the distribution of $\beta_0^\top X$.

We believe that our results are valuable for at least two reasons: (i) introducing a structured nonparametric deconvolution problem where the parametric rate of convergence is attainable, which is uncommon in classical deconvolution problems, and (ii) being able to achieve such rate of convergence without requiring restrictive assumptions on the noise distribution. Another important contribution we make in this work is to propose nonparametric estimators for the unconditional density of $Z$ and for the conditional density of $Z | Y=y$ for any observed response value $y$. Such estimators can be employed to make inference about the latent linear predictor linked to this specific realisation of the response variable.

\subsection{The setting}\label{sec: setting}
Let $X$ be a random vector in $\RR^d$ for $d\geq 1$. Consider $Y$ to be a random variable defined on the same probability space as $X$ such that 
\begin{align}\label{eq: model}
Y = \beta_0^\top X+\epsilon
\end{align}	
where $\epsilon$ is the noise random variable with known distribution and independent of $X$, and $\beta_0\in\RR^d$ is a deterministic vector of unknown coefficients. In this work, we assume that we only have access to $\mathcal X_n=\{X_1,\dots,X_n\}$, a set of independent and identically distributed (i.i.d.)\ copies of $X$ and $\mathcal Y_m = \{Y_1,...,Y_m\}$, a set of i.i.d.\ copies of $Y$. Here, the data $\mathcal X_n$ and $\mathcal Y_m$ can potentially be collected from different sources (i.e., they can be independent),  or originally measured for the same individuals but the link between the responses and covariates is missing,  or anything in between (i.e., there is a non-trivial intersection). Note that it is possible that $n \ne m$.  Still, for ease of exposition, we restrict ourselves to the case of samples of the same size $n$.

Let $\mathcal{B}_0$ be the set of all vectors $\beta\in\RR^d$ such that the unlinked linear model in \eqref{eq: model} holds. In other words,
\begin{align*}
\mathcal{B}_0 =\{\beta\in \RR^d,\ \beta^\top  X \stackrel{d}{=} \beta_0^\top  X\}.
\end{align*}
We will not assume that $\beta_0$ is identifiable, which means that we might be in the case where $| \mathcal B_0 | > 1$.
Under the assumption that $(X,\epsilon)$ is independent of $\{\mathcal X_n,\mathcal Y_n\}$, our goal is to predict the distribution $\mu_0$ of the conditional expectation $\mathbb{E}(Y|X) = \beta^\top _0 X, \ \beta_0 \in \mathcal{B}_0$,  of a future response $Y$ taken from \eqref{eq: model} given  observations $\{\mathcal X_n,\mathcal Y_n\}$ and using our knowledge of the noise distribution $\mu_\epsilon$. We emphasize that the two samples $\mathcal X_n$ and $\mathcal Y_n$ are not linked. We use the word {\it future} to stress the fact that $Y$ is not observed yet.
To do that, we consider $\widehat Z_n:=\widehat\beta_n^\top X^*$ where $\widehat\beta_n$ is a deconvolution least squares estimator as defined in \cite{azadkia2022linear}, and $X^*$ is built as follows: Define $J_n$ to be a random variable which is independent of $\{\mathcal X_n,\mathcal Y_n\}$ and uniformly distributed on $\{1,\dots,n\}$ and let $X^*=X_{J_n}$. Note that $X^*$ is drawn from the empirical distribution of the observations in sample $\mathcal X_n$. This also means that $\widehat Z_n$ is a random draw from the empirical distribution of $\widehat \beta_n^\top X_1, \ldots, \widehat \beta_n^\top X_n$. We consider the distribution of $\widehat Z_n$ as an estimator of the distribution $\mu_0$.

\subsection{Notation}\label{sec: notation} 

All the random variables involved are defined on the same probability space $(\Omega,\mathcal A,\mathbb{P})$ that is equipped with a probability measure denoted by $\mathbb{P}$. We denote the corresponding expectation by $\mathbb{E}$. To alleviate notation, we denote by $P$ and $E$ the conditional probability and expectation given $\{\mathcal{X}_n,\mathcal{Y}_n\}$. Moreover, we will allow ourselves to write $\mathbb{P}(A)$ or $\mathbb{E}(X)$ without justifying that the set $A$ belongs to the sigma-algebra $\mathcal A$ or that $X$ is a measurable function: In case $A\subset\Omega$ is not in $\mathcal{A}$ or the function $X$ is not measurable, then $\mathbb{P}$ and $\mathbb{E}$ must be interpreted as outer probability and outer integral as defined in \cite{vdvwellner23}.

For a given $\beta\in\RR^d$ let 
\begin{align*}
d(\beta,\mathcal B_0):=\inf_{\beta_0\in\mathcal B_0}\|\beta-\beta_0\|
\end{align*}
denote the Euclidean distance of $\beta$ to $\mathcal B_0$.
Moreover, for a given $\beta\in\RR^d$, we will use the following notation:
\begin{itemize}
\item $\mu_\beta^*$ denotes the conditional distribution of the random variable $\beta^\top X^*$ given $\{\mathcal X_n,\mathcal Y_n\}$ (recall it is equal to the empirical distribution of $\beta^\top X_1,\dots,\beta^\top X_n$),
\item $\mu_{\beta}$ denotes the distribution of $\beta^\top X$.
\end{itemize}
Hence, the distribution of $\mathbb{E}(Y|X) = \beta_0^\top X$ from \eqref{eq: model} can be re-written as $\mu_0$ if $\beta_0\in\mathcal B_0$ to emphasize that the distributions $\mu_{\beta_0}$ are all the same in this case. In contrast, note that $\mu_{\beta_0}^*$ might depend on $\beta_0$ even if $\beta_0\in\mathcal B_0$.

Let $\widehat\beta_n$ denote a deconvolution least squares estimator (DLSE) as defined in \cite{azadkia2022linear}, that is 
\begin{align*}
\widehat{\beta}_n \,\, {\in } \,\,\arg \min_{\beta\in\RR^d} \DD_{n}(\beta).
\end{align*}
where
\begin{align*}
\DD_{n}(\beta) := \int\Big( F^Y_{n}(y)-\frac{1}{n}\sum_{i=1}^{n} F^\epsilon(y - \beta^\top  X_i)\Big)^2 d F^Y_{n}(y),
\end{align*}
for $\beta\in\RR^d$. Here,  $F_n^Y$ is the empirical distribution of the sample $\mathcal Y_n$ given by
\begin{align*}
F^Y_{n}(y) := \frac{1}{n} \sum_{i=1}^{n} \one_{(-\infty, y]}(Y_i)
\end{align*}
for all $y\in\RR$, and $F_\epsilon$ is the (known) distribution function of $\epsilon$. 
Using the same notation as above, $\mu^*_{\widehat\beta_n}$ is reserved for the conditional distribution of 
$\widehat\beta_ n^\top X^*$ given $\{\mathcal X_n,\mathcal Y_n\}$, while $\mu_{\widehat\beta_n}$ denotes the conditional  distribution of $\widehat\beta_n^\top X$ given the same samples. Note that randomness in $\mu^*_{\widehat\beta_n}$ and $\mu_{\widehat\beta_n}$ is carried by the uniformly distributed random variable $J_n \in \{1, \cdots, n\}$  and $X$ respectively.

For arbitrary probability measures $\mu$ and $\nu$ on $\RR$ with finite first moment, we denote by $W_1(\mu,\nu)$ the  Wasserstein distance of order 1 between the two measures, that is
\begin{align*}
W_1(\mu,\nu)=\inf_{Z_1\sim\mu,\ Z_2\sim\nu }\mathbb {E}|Z_1-Z_2|  = \int_\RR|F_1(x)-F_2(x)|dx
\end{align*}
with $F_i$ being the distribution function of $Z_i$ {for $i =1,2$}; see Definition 2.1 and Theorem 2.9 in \cite{bobkov2019one}.
We will write $\mu * K$ for the convolution of a measure $\mu$ and a kernel $K$ on $\mathbb{R}$; i.e., 
$(\mu * K)(z) = \int_{\mathbb{R}} K(z - x) d\mu(x).$
Finally, for $M, a>0$ we shall consider the following set of probability measures $\mu$ on $\RR$:  
\begin{align*}
\mathcal C(M,a):=\Big\{\mu:\ \int |x|^{a+2}d\mu(x)\leq M\Big\}.
\end{align*}

\subsection{Organization of the paper}
The paper is organized as follows: Section \ref{sec: main-results} presents the primary contributions of this work. Specifically, we introduce an estimator for the distribution $\mu_0$ and establish its parametric rate of convergence in the Wasserstein distance of order 1. Furthermore, we introduce nonparametric estimators for the unconditional density of $Z$ and the
conditional density of $Z$ given an observed response. This allows us to introduce and study the problem of estimating the true realization of the linear predictor and show how our estimators and theoretical results can be employed in such a statistical task.
Section \ref{sec: num} illustrates several numerical experiments where we show the practical relevance and accuracy of the theoretical findings in various simulation scenarios, providing thereby a clear empirical evidence of their effectiveness.  Section \ref{sec: proofs} provides proofs for the obtained results in which the main arguments leverage concepts and tools from empirical process theory.

\section{Main results} 
\label{sec: main-results}
\subsection{Distributional deconvolution in unlinked linear models}
Consider the statistic
\begin{align*}
\widehat Z_n=\widehat\beta_n^\top X^*
\end{align*}
where $X^*$ follows the empirical distribution based on the sample $\mathcal X_n$  and $\widehat\beta_n$ is a deconvolution least squares estimator (DLSE) as defined in \cite{azadkia2022linear}, see also Section \ref{sec: notation}.
Our estimator for the distribution $\mu_0$ of $\mathbb{E}(Y|X)$ is $\mu_{\widehat\beta_n}^*$, the conditional distribution of $\widehat Z_n$ given observations $\{\mathcal X_n,\mathcal Y_n\}$.
 In other words, $\mu_{\widehat\beta_n}^*$ is the empirical distribution of $\widehat \beta_n^\top X_1, \ldots, \widehat \beta_n^\top X_n$.

It follows from Proposition 1 in \cite{azadkia2022linear} that $\widehat\beta_n$ exists in the sense that (at least) a minimizer of $\DD_{n}$ exists {with probability 1} provided that $n/d > c$ for some fixed constant $c > 0$ and  $F_\epsilon$ is continuous, an assumption that we will also make below. The assumption $n/d > c$ will be automatically satisfied for sufficiently large $n$ since $d$ is assumed to be fixed. Note that $\widehat\beta_n$  might not be unique in which case it will denote any such DLSE, that is any of the elements in $\mathcal{B}_{n}$ where 
\begin{align*}
 \mathcal{B}_{n}:= \{\beta\in\RR^d\text{ s.t. }\DD_{n}(\beta) = \min\limits_{\beta\in\RR^d}\DD_{n}(\beta)\}.
\end{align*}

  We show below that $\mu_{\widehat\beta_n}^*$ is consistent in the sense that $W_1(\mu_{\widehat\beta_n}^*, \mu_0)$ converges to zero in probability. In fact, it converges at the parametric rate $n^{-1/2}$ provided that some appropriate assumptions hold.

A necessary condition for the consistency of any estimator of $\mu_0$ is identifiability of the latter.  
  In our context, with $Z_0  \sim  \mu_0$ independent of $\epsilon$, identifiability of $\mu_0$ means that we have the equivalence
\begin{align}\label{eq: identifiability}
\beta^\top  X +\epsilon\stackrel{d}{=} Z_0+\epsilon \ \ \Longleftrightarrow \ \ \beta^\top  X  \sim \mu_0 
\end{align}
for $\beta\in\RR^d$. 
Although we would like to state our results in great generality and hence only assume that the equivalence in \eqref{eq: identifiability} holds true for all $\beta\in\RR^d$, we give in the following lemma a sufficient condition for (\ref{eq: identifiability}) to be satisfied.

\begin{lemma}\label{lem: identifiability}
Let $\kappa_\epsilon$ denote the characteristic function of the noise $\epsilon$; i.e., 
\begin{align*}
\kappa_\epsilon(t)=\EE\exp(it\epsilon), \  \ t \in \RR.
\end{align*}
If the set of zeros of $\kappa_\epsilon$ does not contain any open and non-empty interval, then the equivalence in \eqref{eq: identifiability} holds  (and whence $\mu_0$ is identifiable).
\end{lemma}

The following theorems extend the results in \cite{azadkia2022linear} about the convergence of $\widehat\beta_n$ to $\beta_0$ when $\beta_0$ is unique to the case where this parameter is possibly non-unique. Also, the theorems establish convergence of the estimator of $\mu_0$ and give sufficient conditions for this convergence to hold at the parametric rate.

\begin{theorem}\label{theo: conv2}
Assume that $\epsilon$ admits a continuous distribution and  that $\PP(\alpha^\top  X = 0) = 0$ for any $\alpha \in \RR^d \setminus \{0\}$. Then, the following holds true.
\begin{enumerate}
    \item  For all $\widehat{\beta}_n\in\mathcal{B}_{n}$ we have that
\begin{align*}
d(\widehat\beta_n,\mathcal B_0)=o_\PP(1),\mbox{ as } n\to\infty.
\end{align*}
    \item If in addition $\mu_0\in \mathcal C(M,a)$ for some $M > 0, a>0$, $\EE\|X\|<\infty$ and the equivalence in \eqref{eq: identifiability} holds for all $\beta\in\RR^d$, then 
\begin{align*}
W_1(\mu_{\widehat\beta_n}^*,\mu_0)=o_\PP(1),\mbox{ as } n\to\infty.
\end{align*}
\end{enumerate}
\end{theorem}

\begin{theorem}\label{theo: rateCV}
Suppose that the distribution function of $\epsilon$ is twice continuously differentiable with bounded first and second derivatives. 
Moreover, suppose that $\EE\|X\|^2<\infty$ and $\PP(\alpha^\top  X = 0) = 0$ for any $\alpha\in \RR^d \setminus \{0\}$. 
For all $\beta\in\RR^d$, define
\begin{align*}
U(\beta)=\int\Big(\int xf^\epsilon(y-\beta^\top x)dF^X(x)\Big)\Big(\int x^\top f^\epsilon(y-\beta^\top x)dF^X(x)\Big)dF^Y(y),
\end{align*}
denote by $\lambda(\beta)$ the smallest eigenvalue of $U(\beta)$, and assume  $\inf_{\beta_0\in\mathcal B_0}\lambda(\beta_0)>0$. Then, the following holds true.
\begin{enumerate}
\item For all $\widehat{\beta}_n\in\mathcal{B}_{n}$, we have that
\begin{align*}
d(\widehat\beta_n,\mathcal B_0)=O_\PP(n^{-1/2}),\mbox{ as } n\to\infty.
\end{align*}
\item If in addition  $\EE\|X\|^p<\infty$ for some $p>2$,  there exists $A>0$ such that for all intervals $[a,b]$ in $\RR$, $\mu_0([a,b])\leq A|b-a|$, and the equivalence in \eqref{eq: identifiability} holds for all $\beta\in\RR^d$, then
\begin{align*}
W_1(\mu^*_{\widehat\beta_n},\mu_0)=O_\PP(n^{-1/2}) ,\mbox{ as } n\to\infty.
\end{align*}
\end{enumerate}
\end{theorem}

Let $\lambda_0:=\inf_{\beta_0\in\mathcal B_0}\lambda(\beta_0)$. The assumption $\lambda_0>0$ is a well-separability condition. In M-estimation, such a condition is rather classical and made to ensure convergence of the M-estimator to the true parameter in case of identifiability. Since the set $\mathcal B_0$ might contain more than one element, the condition $\lambda_0 > 0$ is a natural well-separability requirement in this case as it ensures that $\mathcal D(\beta) > 0$ whenever $\beta \notin \mathcal B_0$, where $\mathcal D$ is the population criterion; i.e.,
\begin{align*}
\mathcal D(\beta):= \int \left(F^Y(y)  -  \int F^\epsilon(y -\beta^\top x) dF^X(x) \right)^2 dF^Y(y).
\end{align*}
The major role played by $\lambda_0$ can be seen in the proof of Theorem \ref{theo: rateCV}: In the first and second claims of the theorem, the $O_\PP(n^{-1/2})$ takes in fact the form $\lambda_0^{-1/2}O_\PP(n^{-1/2})$, respectively $(1+\lambda_0^{-1/2})O_\PP(n^{-1/2})$,  where now, the $O_\PP$-term does not depend on $\lambda_0$. Hence, both bounds on $d(\widehat\beta_n,\mathcal B_0)$ and $W_1(\mu^*_{\widehat\beta_n},\mu_0)$ are of order $n^{-1/2}$ but (for fixed $n$) become larger for smaller values of $\lambda_0$.

Note that the assumption that $\mu_0\in\mathcal C(M,a)$  stated in Theorem \ref{theo: conv2} is not needed in Theorem \ref{theo: rateCV}  because in that theorem we require the stronger condition $\EE\|X\|^p<\infty$ for some $p>2$. Indeed, note that  $\mu_0\in\mathcal C(M,a)$ is equivalent to $\EE|\beta_0^\top X|^{a+2}\leq M$ for all $\beta_0\in\mathcal{B}_0$. Using the Cauchy-Schwarz inequality implies that
\begin{align*}
\EE|\beta_0^\top X|^{a+2}\leq \|\beta_0\|^{a+2}\EE\|X\|^{a+2}.
\end{align*}
Thanks to Lemma \ref{lem: projection} in Section \ref{sec: proofs} below, there exists $B>0$ such that $\|\beta_0\|\leq B$ for all $\beta_0\in\mathcal B_0$, thus, taking $a:=p-2>0$ and $M:=B^{a+2}\EE\|X\|^p<\infty$, we get $\EE|\beta_0^\top X|^{a+2}\leq M$.
In other words, the assumption $\EE\|X\|^p<\infty$ for some $p>2$ implies that we can find $a>0$ and $M>0$ such that $\mu_0\in\mathcal C(M,a)$.
 
\subsection{Density deconvolution and inference on the latent predictor}
A relevant inferential task involves the scenario where we observe a realization $y_0$ from the response variable $Y$ such that  $Y = \beta_0^\top  X + \epsilon$. Suppose that the realization of the covariate $X$, which generated $y_0$,  is non-observed. At the same time, suppose that we are provided with a sample of unlinked covariates and responses $\{\mathcal{X}_n, \mathcal{Y}_n\}$. Let $Z := \mathbb{E}(Y|X) = \beta_0^\top  X$ and $z_0$ be the true value of $Z$ which is linked to $y_0$ but not accessible. Then, it is of interest to make inference about $z_0$ conditionally on the event $\{Y = y_0\}$. One possible approach is based on estimating the conditional density of $Z|Y=y_0$, together with several summaries of the corresponding conditional distribution, including location measures such as the mean and mode, as well as quantiles.

In our previous results, we introduced an estimator for the distribution of $Z$ based on the empirical distribution of $\widehat{\beta}_n^\top  X_1, \ldots, \widehat{\beta}_n^\top  X_n$. Assuming that the distribution of $Z$ admits a density with respect to Lebesgue measure, $f_{Z}$ say, we can apply a Kernel Density Estimator (KDE) with a given kernel and bandwidth and which is based on $\widehat{\beta}_n^\top  X_i, i =1, \ldots, n$. We denote the resulting density estimator by $\widehat{f}_Z$. 
Let $f_Y$, $f_{Y|Z}$ and $f_{Z|Y}$ denote the density of $Y$, and conditional densities of $Y$ given $Z$ and $Z$ given $Y$, respectively. Using the Bayes rule, we have that
\begin{align}
\label{eq: fZ|Y}
f_{Z|Y}(z|y) = \frac{f_{Y|Z}(y|z) f_Z(z)}{f_Y(y)}, \ \  y, z \in \mathbb R.
\end{align}
Since $f_{Y|Z}(y|z)  = f_\epsilon(y - z)$, where $f_\epsilon$ is the known density of the noise $\epsilon$, the only missing element is an estimator for the (marginal) density $f_Y$, for which a simple estimator is given by
\begin{align}
	\label{eq: f_Y_first}
	\widehat{f}_Y(y) = \frac{1}{n} \sum_{i=1}^n f_\epsilon(y - \widehat{\beta}_n^\top  X_i), \ \ y \in \mathbb R.
\end{align}
Thus, we consider the following estimator of the conditional density $f_{Z|Y}$:
\begin{align}
	\label{eq: estimator-Z|Y}
	\widehat{f}_{Z|Y}(z|y) = \frac{f_\epsilon(y - z) \widehat{f}_Z(z)}{\widehat{f}_Y(y)}, \ \  y, z \in \mathbb R.
\end{align}

In the following theorem we provide guarantees for both the unconditional and conditional density estimators.
\begin{theorem}
\label{theo: kde}
     Suppose $\mu_0$ admits a density $f_Z\in C^2(\mathbb{R})$ with bounded first and second derivatives.
     Moreover, for some $L > 0$ consider an $L$-Lipschitz kernel $K:\mathbb{R}\to [0, \infty)$ such that $\int K(u) du=1$, $\int uK(u) du=0$ and $\int u^2 K(u) du<\infty$. Let $K_h(z):= h^{-1}K(z/h)$ with $h=h_n\asymp n^{-1/8}$. Under the assumptions of Theorem \ref{theo: rateCV}, the following holds true.
    \begin{enumerate}
        \item The estimator $\widehat{f}_{Z}(z):=(\mu^*_{\widehat\beta_n} * K_h) (z)=  n^{-1} \sum_{i=1}^n  K_h(z- \widehat\beta_n^\top X_i)$
        satisfies 
        $$
        \|\widehat{f}_{Z}-f_Z\|_{\infty} = O_\mathbb{P}(n^{-1/4}).
        $$

        \item If, moreover, $f_\epsilon$ is bounded and $M$-Lipschitz for some $M > 0$, and  $f_Y(y_0) \in (0, \infty)$, then  
        the estimator $\widehat{f}_{Z|Y}(z|y_0) = f_\epsilon(y_0 - z) \widehat{f}_Z(z)/\widehat{f}_Y(y_0)$, with $\widehat f_Y$ given in \eqref{eq: f_Y_first}, satisfies
        $$
        \|\widehat{f}_{Z|Y}(\cdot|y_0)-f_{Z|Y}(\cdot|y_0)\|_{\infty} = O_\mathbb{P}(n^{-1/4}).
        $$
    \end{enumerate}
\end{theorem}
It is worth remarking that in the setting of Theorem \ref{theo: kde}, under the additional assumptions $f_Z\in C^s(\RR)$, $\int u^j K(u) du = 0$ for $j=1\dots s-1$, $\int u^s K(u) du <\infty$, and $h_n\asymp n^{-1/(2s+4)}$, it is possible to obtain  the rate of convergence $n^{-s/(2s+4)}$ in sup-norm for both density estimators $\widehat f_Z$ and $\widehat{f}_{Z|Y}(\cdot|y_0)$. Moreover, it is evident that the resulting rate of convergence is slower than the one achieved in classical KDE theory. The question of optimality lies beyond the scope of this work and is left for future research.

Now suppose that $Y = y_0$, a realization from the distribution of $\beta_0^T X + \epsilon =  z_0 + \epsilon$. The true mean $z_0$ is unknown and the goal is to use the estimator of the conditional density given in \eqref{eq: estimator-Z|Y} to make inference about $z_0$ given that $Y=y_0$ was observed. We consider
\begin{align}
	\label{eq: point-estimators-Z|Y}
	\widehat{E}(Z|Y=y_0) = \int z \, \widehat{f}_{Z|Y}(z|y_0) \, dz, \quad \quad \widehat{M}(Z|Y=y_0) = \arg \max_{z} \, f_\epsilon(y_0 - z)\widehat{f}_Z(z),
\end{align}
which can be computed via numerical methods.  In a Bayesian modelling framework, interval estimators can also be considered. These are derived as credible intervals based on the quantiles of the distribution of $Z|Y=y_0$, which are in turn obtained via $\widehat f_{Z|Y=y_0}$. {Let $Q_{Z|Y}(\cdot |y_0)$  be the quantile function of $Z|Y=y_0$. For a given probability $p\in (0,1)$, consider the estimator
    \begin{align}\label{eq: quantile}
    \widehat Q_{Z|Y}(p|y_0):= \inf \Big\{x\in\mathbb{R}:\int_{-\infty}^x \widehat f_{Z|Y}(z|y_0)dz \ge p \Big\}  =  \widehat F_{Z|Y}^{-1}(p|y_0). 
    \end{align}
where $\widehat F_{Z|Y}(\cdot|y_0) = \int_{-\infty}^{\cdot} \widehat f_{Z|Y}(z|y_0) dz$. 
The following Proposition provides rates of convergence for \eqref{eq: point-estimators-Z|Y} and \eqref{eq: quantile}.
\begin{proposition}
\label{prop: mean-mode-quantile}
    Suppose that the assumptions of Theorem \ref{theo: kde} are satisfied and let $y_0\in \RR$ such that $f_Y(y_0) \in (0, \infty)$.
    
    \begin{enumerate}
        \item For the mean estimator it holds 
        $$|\widehat{E}(Z|Y=y_0)-\EE(Z|Y=y_0)| = O_\PP(n^{-1/4}). $$
        \item  Suppose that $z \mapsto f_{Z|Y}(z| y_0)$ admits a maximizer at $z^*:=\MM(Z|Y=y_0)$ such that if $\vert z  -  z^* \vert > M $  then for some $\delta > 0$ (which depends on $M$) it holds that
        $$
        f_{Z|Y}(z^*|y_0)-f_{Z|Y}(z|y_0)\geq \delta.
        $$
        Then, $\widehat M(Z|Y=y_0)  \overset{\PP}{\to}  z^*$.  Furthermore, if 
         $$f_{Z|Y}(z^*|y_0)-f_{Z|Y}(z|y_0)\geq \lambda|z^*-z|^\alpha,$$
         on a small neighbourhood of $z^*$ for some constants $\alpha, \lambda>0$, then
        $$|\widehat M(Z|Y=y_0) - \MM(Z|Y=y_0)| = O_\PP(n^{-1/(4\alpha)}).$$
        \item Let $p\in (0,1)$ and suppose $f_{Z|Y}(\cdot|y_0)$ is strictly positive in a neighbourhood of $Q_{Z|Y}(p|y_0)$. Moreover, assume that there exists constants $c_1,c_2,c_3>0$ such that $f_{Z|Y}(z|y_0)\leq c_2e^{-c_3|z|}$ for all $|z|>c_1$. Then, 
        $$|\widehat Q_{Z|Y}(p|y_0)-Q_{Z|Y}(p|y_0)|= O_{\mathbb{P}}(n^{-1/4}\log n).$$
              
    \end{enumerate}
\end{proposition}

In the next section, we show numerically that when the noise $\epsilon$ has a small variance $\sigma^2$, the density estimator $\widehat{f}_{Z|Y=y_0}$ and the estimators in \eqref{eq: point-estimators-Z|Y} perform well in recovering the true value $z_0$ linked to $y_0$.

\section{Numerical experiments}\label{sec: num}
In this section, we illustrate the theoretical results using simulated data.  Our first goal is to show that $\mu_{\widehat \beta_n}^*$ is close to $\mu_0$ in the sense of $W_1$-distance and that the estimation error decreases at a parametric rate. In addition, we aim to evaluate the performance of the introduced estimators both qualitatively and quantitatively. For all experiments, the estimator $\widehat{\beta}_n$ is calculated using the $\texttt{optim}$ function from the R-package $\texttt{stats}$. The simulation study was conducted using 500 independent Monte Carlo replications, and the results were reported in the form of Monte Carlo averages of the relevant quantities.
Let $\mathcal{N}(\mu,\Sigma)$ denote the multivariate Gaussian distribution with mean $\mu \in \mathbb R^d$ and positive definite covariance matrix $\Sigma \in \mathbb{R}^{d \times d}$.  In addition,  we denote by $\mathrm{Gamma}(a,b)$ the Gamma distribution with shape $a >  0$ and scale $b > 0$ respectively.
In all simulations we took $\epsilon \sim \mathcal{N}(0,\sigma^2)$, for $\sigma > 0$, and considered four different settings:
\begin{enumerate}[label=(\alph*)] 
	\item $d=2$, $X \sim \mathcal{N}(0,I)$, $\beta_0^\top  = (3,-5) $,
		\item $d=3$, $X \sim \mathcal{N}(0,I)$, $\beta_0^\top  = (-1.5,2,7) $,
	\item $d=2$, $X=
	(X^{(1)},X^{(2)})^\top $ with $X^{(1)} \sim \mathrm{Gamma}(1,1)$, $X^{(2)}\sim \mathrm{Gamma}(2,4)$,  $\beta_0^\top  = (1,2) $, 
\item $d=3$, $X=
	(X^{(1)},X^{(2)},X^{(3})^\top $ with $X^{(1)} \sim \mathrm{Gamma}(1,1)$, $X^{(2)}\sim \mathrm{Gamma}(2,4)$, $X^{(3)}\sim \mathrm{Gamma}(1.5,3)$, $\beta_0^\top  = (0.5,2,3) $.
\end{enumerate}

Note that in the simulation settings (a) and (b) it holds that $\mathcal B_0 =  \{ \beta \in \mathbb R^2:  \Vert \beta \Vert = \Vert \beta_0 \Vert =  \sqrt{34}\} $ for (a)  and $\mathcal B_0 =  \{ \beta \in \mathbb R^3:   \Vert \beta \Vert = \Vert \beta_0 \Vert = \sqrt{55.25}   \} $ for (b). Hence, $\beta_0$ is not identifiable in these settings.
For the settings in (c) and (d), we can use Theorem 5 of \cite{fadmartjohn}. This theorem implies that if $X = (X_1, \ldots, X_p)^\top$ such that $X_1, \ldots, X_p$ are independent and $X_i \sim \text{Gamma}(a_i, b_i), i =1, \ldots, p$ with the property that 
\begin{align}\label{Cond}
\sum_{i \in I} a_i = \sum_{j \in J} a_j  \  \Longrightarrow \   I = J
\end{align}
for $I, J \subset \{1, \ldots, p \}$, then $\vert \mathcal B_0 \vert =1$. In  setting (c), we have that $a_1 = 1 \ne  a_2 =2$, and the condition in (\ref{Cond}) is clearly satisfied. In setting (d), we have that $a_1 =1, a_2 = 2, a_3= 1.5$, and we can easily check that (\ref{Cond}) holds.
Thus, in settings (c) and (d), $\beta_0$ is identifiable.

\medskip 
\noindent \textbf{Experiment 1.} We fix $\sigma=1$ and calculate $W_1(\mu_{\widehat \beta_n}^*,\mu_0)$ on the basis of i.i.d.\ samples with increasing size $n \in \{1000,2000,\dots,5000\}$, drawn from the distribution of each of the four settings (a)-(d). 
For this purpose, a general method to practically compute $W_1(\mu_{\widehat \beta_n}^*,\mu_0)$ is the following: First, draw a large sample from $\mu_0$ of size $m=10^6$ and compute the associated empirical distribution function $ F_{0,m}$. Successively, compute the empirical distribution function $F^*_{\widehat \beta_n}(x)=\frac{1}{n}\sum_{i=1}^n \one_{\widehat \beta_nX_i\leq x}$ based on the sample $\widehat \beta_n^\top X_1,\dots, \widehat \beta_n^\top X_n$. Now, introduce the sequence $\{x_j\}_{j=1}^{n+m}$ in $\mathbb{R}$ obtained merging and sorting the two empirical samples and the sequence $\{u_j\}_{j=1}^{n+m+1}=\{0\}\cup  \{v(x_j)\}_{j=1}^{n+m} $ in $[0,1]$ such that 
\begin{align*}
v(x_j) =
\begin{cases}
F^*_{\widehat \beta_n}(x_j), & \text{if } x_j \text{ belongs to the sample from } \mu_{\widehat \beta_n}^*, \\
F_{0,m}(x_j), & \text{if } x_j \text{ belongs to the sample from } \mu_0,
\end{cases}
\end{align*}
namely, $\{u_j\}_{j=1}^{n+m+1}$ contains all the sorted values attained by the two empirical distribution functions.
Finally, consider the following numerical approximation
\begin{align*}
W_1(\mu_{\widehat \beta_n}^*,\mu_0) \approx \sum_{j=1}^{n+m+1}  |F^{*-1}_{\widehat \beta_n}(u_{j+1})-F^{-1}_{0,m}(u_{j+1})|(u_{j+1}-u_j),
\end{align*}
which is arbitrarily accurate by increasing $m$.

Now, denote by $W_1^{(j)}(\mu_{\widehat \beta_n}^*,\mu_0)$ the $j$-th Monte Carlo replication of $W_1(\mu_{\widehat \beta_n}^*,\mu_0)$ for $j =1, \ldots, J=500$. Based on these replications, we compute the slope of the log-linear fit (linear fit on the logarithmic scale) for 
\begin{align*}
\overline W_1^{(k)} := \frac{1}{J}\sum_{j=1}^J \Big(W_1^{(j)}(\mu_{\widehat \beta_n}^*,\mu_0)\Big)^k, \quad k=1,2,3, 
\end{align*}
and also the $99\%$ empirical quantile, which we denote by $q_{W_1}^{0.99}$. The results are reported in Table \ref{tab: slopes} and suggest that $\overline W_1^{(k)} \asymp n ^{-k/2}$ and $q_{W_1}^{0.99} \asymp n ^{-1/2} $, thus providing strong empirical evidence for the theory.
Figure \ref{fig: rate-plot} also reports graphically the Monte Carlo average $\overline W_1^{(1)}$ of the Wasserstein distances on the logarithmic scale.

\begin{table}[!ht]
\centering
\caption{Estimated slopes of the log-linear fit for the first 3 moments and 99\% quantile of the Monte Carlo distribution of $W_1(\mu_{\widehat \beta_n}^*,\mu_0)$ under the settings (a)-(d).}
\begin{tabular}{@{}l|cccc@{}}
\toprule
Setting & $\overline W_1^{(1)}$ & $\overline W_1^{(2)}$ & $\overline W_1^{(3)}$ & $q_{W_1}^{0.99}$  \\ \midrule
(a)     & -0.502 & -1.007 & -1.516 & -0.567 \\
(b)     & -0.493 & -0.987 & -1.479 & -0.514 \\
(c)     & -0.501 & -1.000 & -1.496 & -0.464 \\
(d)     & -0.485 & -0.953 & -1.365 & -0.477 \\ \bottomrule
\end{tabular}
\label{tab: slopes}
\end{table}

\begin{figure}[!ht]
	\centering
	\includegraphics[width=\textwidth]{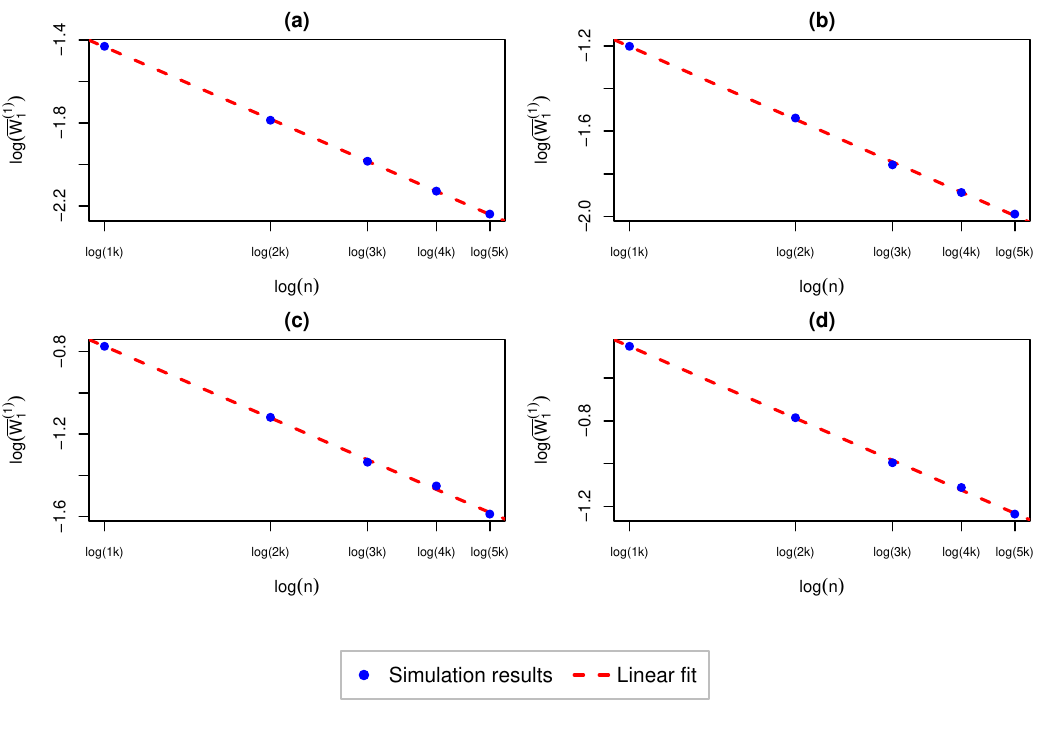}
	\caption{Log-linear fit for the first moment $\overline W_1^{(1)}$ of the Monte Carlo distribution of $W_1(\mu_{\widehat \beta_n}^*,\mu_0)$ under the settings (a)-(d) for varying sample size. The slopes of the linear fit are -0.502 in (a), -0.493 in (b),  -0.501 in (c) and -0.485 in (d).
    }
	\label{fig: rate-plot}
\end{figure}

\medskip

\noindent \textbf{Experiment 2.}  In all the four settings, we again fix $\sigma=1$.  Here, a random sample of size $n=500$ is drawn from each of the corresponding distributions and then used to calculate the KDE $\widehat f_Z$ introduced in Theorem \ref{theo: kde} using a Gaussian Kernel with a bandwidth systematically chosen via the simple rule of thumb $h = n^{-1/8} \operatorname{sd}(\widehat \beta^\top_n X_1,\dots,\widehat \beta^\top_n X_n )$ where $\operatorname{sd}$ denotes the sample standard deviation, in order to ensure that it is appropriately rescaled and matches the setting of Theorem \ref{theo: kde}.  For a given value of $z_0$ randomly generated from $Z$, a new random response $Y \sim  z_0  + \epsilon$ was generated. Based on the realization $y_0$ of $Y$, 
we compute the estimator $\widehat f_{Z|Y=y_0}:=\widehat f_{Z|Y}(\cdot|y_0)$ for the density of the random variable $Z|Y=y_0$ given in \eqref{eq: estimator-Z|Y}, the resulting mean $\widehat{E}(Z|Y=y_0)$ and mode $ \widehat{M}(Z|Y=y_0)$ given in \eqref{eq: point-estimators-Z|Y}, and the resulting 95\% credible interval that we denote by $\mathrm{CI}_{0.95}( \widehat f_{Z|Y=y_0})$.
We remark that when $Z$ has a small variance, adaptive quadrature methods for computing $\widehat{E}(Z|Y=y)$ may perform poorly. For this reason, we resort to importance sampling to obtain a more robust and stable approximation of the integral in \eqref{eq: point-estimators-Z|Y} and of the credible intervals based on $\widehat f_{Z|Y=y}$. Moreover, most density evaluations are performed on the logarithmic scale to enhance numerical stability by reducing errors in floating-point arithmetic.
The results of the experiment are shown in Figure \ref{fig: estimators-plot}. In all settings, it is clear that $\widehat f_Z$ recovers well the true $f_Z$.  However, as expected, the conditional density estimator $\widehat f_{Z|Y=y_0}$ yields a much better prediction for the true $z_0$.

\begin{figure}[!ht]
	\centering
	\includegraphics[width=\textwidth]{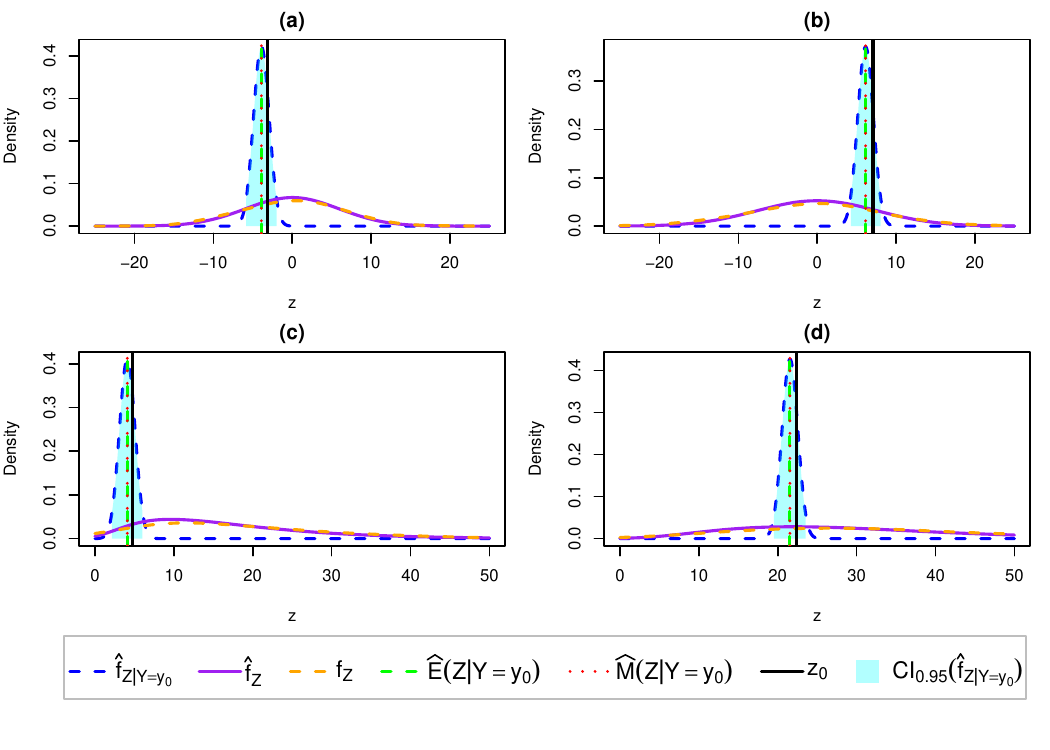}
	\caption{Graphical results of Experiment 2 under the settings (a)-(d). The true $(z_0,y_0)$ are $(-3.155, -3.996)$ in (a), $(7.074, 6.233)$ in (b), $(4.809, 3.968)$ in (c) and $(22.367, 21.526)$ in (d).}
	\label{fig: estimators-plot}
\end{figure}

\medskip
\noindent \textbf{Experiment 3.} 
For the four settings, we fix $\sigma=1$ and let $n\in \{50,100,500\}$. Here, we shall compare the performance at recovering the true $z$ of the point estimators $\widehat{E}(Z|Y=y) $ and $ \widehat{M}(Z|Y=y)$ with the estimated mean and mode of $Z$, respectively given by 
\begin{align*}
\widehat E (Z) = \frac{1}{n}\sum_{i=1}^n \widehat \beta_n^\top X_i, \quad \widehat M(Z) = \arg\max_z \widehat f_Z(z).
\end{align*}
The performance comparison is done by computing the performance ratios 
\begin{align*}
R_{\widehat{E}}=\frac{\mathrm{MSE}(\widehat{E}(Z|Y=y))}{\mathrm{MSE}(\widehat{E}(Z))}, \quad  R_{\widehat{M}}=\frac{\mathrm{MSE}(\widehat{M}(Z|Y=y))}{\mathrm{MSE}(\widehat{M}(Z))}.
\end{align*}
To explain how the MSE is computed, let us focus on the conditional mode. For each Monte Carlo replication,  we compute this mode using a sample of size $n$, thereby obtaining $y \mapsto \widehat{M}^{(j)}(Z|Y=y)$ for $j =1, \ldots J=500$. The corresponding MSE is then obtained as the average squared difference between the mode and the true $z$, 
\begin{align}
\label{eq: MSE}
\mathrm{MSE}(\widehat{M}(Z|Y=y))=\frac{1}{J}\sum_{j=1}^J\frac{1}{T} \sum_{t=1}^T  \Big({ \widehat{M}^{(j)}(Z|Y=y_t)-z_t}\Big)^2
\end{align}
where $T$ is the size of test sample $\{(z_t, y_t)\}_{t=1}^T$,  which is set to be equal to 100. A similar approach is used for the conditional mean and the unconditional estimators.
In this experiment, we also compare the performance of two credible intervals, one obtained through the symmetric extreme empirical quantiles based on $\widehat \beta_n^\top X_1,\dots, \widehat \beta_n^\top X_n$, and the other obtained through importance sampling based on $\widehat f_{Z|Y=y}$. For a given nominal $\alpha \in (0,1)$, with a slight abuse of notation, we denote such intervals by $ \mathrm{CI}_{1-\alpha}(\widehat f_{Z})$ and $\mathrm{CI}_{1-\alpha}(\widehat f_{Z|Y=y})$, respectively.
For the experiment, we set $\alpha=0.05$ and report the Monte Carlo empirical frequentist coverage of the credible intervals computed as
\begin{align*}
    \mathcal{C}(\mathrm{CI}_{1-\alpha}(\widehat f_{Z|Y=y}))= \frac{1}{J}\sum_{j=1}^J\frac{1}{T}\sum_{t=1}^T\one\{z_t\in \mathrm{CI}_{1-\alpha}(\widehat f_{Z|Y=y_t}^{(j)})\},
\end{align*}
and similarly for $ \mathrm{CI}_{1-\alpha}(\widehat f_{Z})$.
Moreover, we report the length ratio 
\begin{align*}
R_{\mathrm{len} (\mathrm{CI}_{1-\alpha})}= \frac{\frac{1}{J}\sum_{j=1}^J\frac{1}{T} \sum_{t=1}^T    \mathrm{len} (\mathrm{CI}_{1-\alpha}(\widehat f^{(j)}_{Z|Y=y_t}))}{\frac{1}{J}\sum_{j=1}^J\frac{1}{T} \sum_{t=1}^T   \mathrm{len}  (\mathrm{CI}_{1-\alpha}(\widehat f^{(j)}_{Z}))},
\end{align*}
where $\mathrm{len} ([a,b])=b-a$.
The results are reported in Table \ref{tab: comparison}. It is at once apparent that the proposed conditional estimators are preferable in terms of $\mathrm{MSE}$. Moreover, for all scenarios, both interval estimators attain the correct frequentist coverage, with the conditional one being calibrated already for very small sample sizes and leading to much narrower intervals.

\begin{table}[!ht]
\centering
\caption{Performance comparison under the settings (a)-(d) for varying sample size.}
\begin{tabular}{@{}ll|ccccc@{}}
\toprule
Setting & $n$ & $R_{\widehat{E}}$ & $R_{\widehat{M}}$ 
& $\mathcal{C}(\mathrm{CI}_{0.95}(\widehat f_{Z|Y=y}))$ 
& $\mathcal{C}(\mathrm{CI}_{0.95}(\widehat f_{Z}))$ 
& $R_{\mathrm{len} (\mathrm{CI}_{1-\alpha})}$ \\ \midrule

(a) &  50  & 0.028 & 0.027 & 0.949 & 0.915 & 0.182 \\
    & 100  & 0.029 & 0.028 & 0.949 & 0.933 & 0.175 \\
    & 500  & 0.029 & 0.029 & 0.949 & 0.946 & 0.170 \\ \midrule

(b) &  50  & 0.018 & 0.017 & 0.950 & 0.912 & 0.144 \\
    & 100  & 0.018 & 0.017 & 0.950 & 0.934 & 0.138 \\
    & 500  & 0.018 & 0.018 & 0.950 & 0.946 & 0.134 \\ \midrule

(c) &  50  & 0.009 & 0.007 & 0.948 & 0.913 & 0.101 \\
    & 100  & 0.008 & 0.007 & 0.948 & 0.931 & 0.096 \\
    & 500  & 0.008 & 0.006 & 0.948 & 0.946 & 0.092 \\ \midrule

(d) &  50  & 0.010 & 0.010 & 0.949 & 0.910 & 0.105 \\
    & 100  & 0.010 & 0.010 & 0.949 & 0.929 & 0.102 \\
    & 500  & 0.010 & 0.010 & 0.948 & 0.946 & 0.098 \\ \bottomrule

\end{tabular}
\label{tab: comparison}
\end{table}

\medskip
\noindent \textbf{Experiment 4.} Under the four settings, for a given linked realization $(z_0, y_0)$ we use $\widehat{E}(Z|Y=y_0) $ and $ \widehat{M}(Z|Y=y_0)$ in \eqref{eq: estimator-Z|Y} and \eqref{eq: point-estimators-Z|Y} for predicting $z_0$, and we evaluate again their predictive performance by estimating their associated Mean Squared Error (MSE) over Monte Carlo replications of $(z, y)$. As opposed to the previous experiment we perform such evaluation for different values of $n$ and also $\sigma$ (recall that $\sigma =1$ in Experiment 3). The MSE of the estimators is computed as in \eqref{eq: MSE}. Here, we take $n \in \{50,100,500\}$ and $\sigma^2 \in \{0.5,1,1.5,2,2.5\}$. 
The results are reported in Tables \ref{tab: mean} and  \ref{tab: mode} where it is shown that the proposed estimators of the true value of the linear predictor are comparable in terms of their performance. In both tables, the MSE decreases as $n$ increases under all the considered settings. However, it does not tend to $0$ since the magnitude of the noise constraints the accuracy that any estimator can achieve. The MSE is comparable to $\sigma^2$, and it is interesting to note that noise introduces a form of aleatoric uncertainty that cannot be reduced as opposed to epistemic uncertainty that decreases with increasing sample size; see \cite{hullermeier2021aleatoric}.

\begin{table}[!ht]
\centering
\caption{MSE for the mean estimator under the settings (a)-(d) for varying $\sigma^2$ and sample size.}
\begin{tabular}{@{}ll|ccccc@{}}
\toprule
Setting & $n$ & $\sigma^2=0.5$ & $\sigma^2=1$ & $\sigma^2=1.5$ & $\sigma^2=2$ & $\sigma^2=2.5$ \\ \midrule

(a) &  50  & 0.497 & 0.981 & 1.454 & 1.915 & 2.365 \\
    & 100  & 0.496 & 0.980 & 1.452 & 1.911 & 2.359 \\
    & 500  & 0.497 & 0.979 & 1.449 & 1.906 & 2.352 \\ \midrule

(b) &  50  & 0.499 & 0.989 & 1.472 & 1.948 & 2.416 \\
    & 100  & 0.498 & 0.988 & 1.469 & 1.942 & 2.408 \\
    & 500  & 0.498 & 0.987 & 1.468 & 1.940 & 2.405 \\ \midrule

(c) &  50  & 0.645 & 1.138 & 1.627 & 2.111 & 2.590 \\
    & 100  & 0.551 & 1.044 & 1.532 & 2.015 & 2.493 \\
    & 500  & 0.517 & 1.010 & 1.496 & 1.976 & 2.449 \\ \midrule

(d) &  50  & 0.546 & 1.041 & 1.532 & 2.017 & 2.499 \\
    & 100  & 0.512 & 1.006 & 1.496 & 1.980 & 2.460 \\
    & 500  & 0.509 & 1.002 & 1.491 & 1.974 & 2.453 \\ \bottomrule

\end{tabular}
\label{tab: mean}
\end{table}

\begin{table}[!ht]
\centering
\caption{MSE for the mode estimator under the settings (a)-(d) for varying $\sigma^2$ and sample size.}
\begin{tabular}{@{}ll|ccccc@{}}
\toprule
Setting & $n$ & $\sigma^2=0.5$ & $\sigma^2=1$ & $\sigma^2=1.5$ & $\sigma^2=2$ & $\sigma^2=2.5$ \\ \midrule

(a) &  50  & 0.496 & 0.980 & 1.453 & 1.914 & 2.364 \\
    & 100  & 0.496 & 0.979 & 1.451 & 1.910 & 2.358 \\
    & 500  & 0.496 & 0.978 & 1.447 & 1.905 & 2.350 \\ \midrule

(b) &  50  & 0.498 & 0.988 & 1.471 & 1.947 & 2.415 \\
    & 100  & 0.497 & 0.987 & 1.467 & 1.941 & 2.407 \\
    & 500  & 0.497 & 0.986 & 1.466 & 1.938 & 2.402 \\ \midrule

(c) &  50  & 0.643 & 1.134 & 1.622 & 2.105 & 2.584 \\
    & 100  & 0.550 & 1.042 & 1.529 & 2.011 & 2.488 \\
    & 500  & 0.515 & 1.007 & 1.492 & 1.972 & 2.446 \\ \midrule

(d) &  50  & 0.545 & 1.039 & 1.530 & 2.015 & 2.497 \\
    & 100  & 0.510 & 1.004 & 1.493 & 1.978 & 2.458 \\
    & 500  & 0.507 & 1.001 & 1.489 & 1.971 & 2.449 \\ \bottomrule

\end{tabular}
\label{tab: mode}
\end{table}

\section{Proofs} \label{sec: proofs}

\subsection{Preparatory lemmas}

The following lemma is an extension of Theorem 3.2.5 in \cite{vdvwellner23}.  
\begin{lemma}\label{lem: rateofCV}
Let $\Theta$ be equipped with a semi-metric $d$ with $\mathcal{B} \subseteq \Theta$ and define $d(\theta, \mathcal{B})  =  \inf_{\theta' \in \mathcal{B}} d(\theta, \theta')$.
Let $\MM_n$ be a stochastic process indexed by $\Theta$ and $\MM:\Theta\mapsto\RR$ a deterministic function such that for every  $\theta_0\in\cal B$ and $\theta\in\Theta$ satisfies
\begin{align*}
\MM(\theta)-\MM(\theta_0) \leq -C_0^2d^2(\theta,\mathcal B)
\end{align*}
for some $C_0>0$.
Suppose that for all $n$ and sufficiently small $\delta   >  0$, the centred process $\MM_n-\MM$ satisfies
\begin{align}\label{eq: rateofCV}
\mathbb{E}\sup_{(\theta,\theta_0) \in \Theta \times \mathcal{B} :\ d(\theta,\theta_0)<\delta}\Big|(\MM_n-\MM)(\theta)-(\MM_n-\MM)(\theta_0)\Big|\lesssim\frac{\phi_n(\delta)}{\sqrt n},
\end{align}
for functions $\phi_n$ such that $\delta\mapsto\phi_n(\delta)/\delta^\alpha$ is decreasing for some $\alpha<2$ (not depending on $n$). Moreover, suppose that for all $n$,
\begin{align*}
r_n^2\phi_n\Big(\frac1{r_n}\Big)\leq\sqrt n,
\end{align*}
the sequence $\widehat\theta_n$ satisfies 
\begin{align*}    
\MM_n(\widehat\theta_n)\geq\sup_{\theta_0\in\cal B}\MM_n(\theta_0)-O_\PP(r_n^{-2}),
\end{align*}
and $d(\widehat\theta_n,\mathcal B)\to 0$ in outer probability, then 
\begin{align*}
    r_nd(\widehat\theta_n,\mathcal B)= C_0^{-1}O_{\PP}(1)
\end{align*}
where the $O_{\PP}$-term does not depend on $C_0$. 
\end{lemma}
In the next lemma, we state the compactness of the set $\mathcal B_0$.
\begin{lemma}\label{lem: projection}
Suppose that $E\| X\|^2 < \infty$. Then, $\mathcal B_0$ is compact and  for every $\beta\in\RR^d$ there exists $\beta_0\in\mathcal B_0$ such that 
\begin{align*}
\inf_{\beta'\in\mathcal{B}_0}\|\beta-\beta'\|=\|\beta-\beta_0\|.
\end{align*}
\end{lemma}

The next lemma is a slight extension of Proposition 11 in \cite{azadkia2022linear} to the case where $| \mathcal B_0 | > 1$, and we use the empirical process notation \citep{vdvwellner23}. With $V$ denoting either $X$ or $Y$, we introduce
\begin{align*}
\GG_n^V:=\sqrt n(\PP_n^V-\PP^V)
\end{align*}
where $\PP^V$ denotes the distribution of the random variable $V$ and $\PP_n$ denotes the empirical distribution of $n$ independent random variables distributed as $V$, and for a class of functions $\mathcal F$ we let
\begin{align*}
\|\GG_n^V\|_{\mathcal F}:=\sup_{f\in\mathcal F}|\GG^V_n f|.
\end{align*}
Moreover, if $\mathcal F$ has envelop $F$ we define
\begin{align*}
J(\delta,\mathcal F):=\sup_{\QQ}\int_0^\delta\sqrt{1+\log N(\eta\|F\|_\QQ,\mathcal F,L^2(\QQ))}d\eta
\end{align*}
where $N(\nu,\mathcal F,\|\cdot\|)$ is the $\nu$-covering number of $\mathcal F$ with respect to $\|\cdot\|$, and where the supremum is taken over  all probability measures $\QQ$ such that $\|F\|_\QQ^2:=\int F^2d\QQ>0.$

\begin{lemma}\label{lem: entropy}
Suppose that $\mathbb E\| X\|^2 < \infty$, and that $\epsilon$ has a bounded density.
Consider the classes of functions $\mathcal{M}$ and $\mathcal{G}_\delta$, where $\delta>0$, defined as
\begin{align*}
\mathcal{M}:= \{ m_{\beta_0, y},\ y \in \RR, \beta_0\in\mathcal B_0  \},
\end{align*}
with $m_{\beta, y}(x) = F^\epsilon(y - \beta^\top x), \ \beta \in \RR^d$, and 
\begin{align*}
\mathcal{G}_\delta := \{g_{\beta,\beta', y}(x) =  m_{\beta, y}(x) - m_{\beta', y}(x),\ \|\beta - \beta'\|\le\delta, y\in\RR \}
\end{align*}
Then, we have
\begin{align*}
\EE\| \mathbb{G}^X_n\|^2_{\mathcal{M}}  \lesssim 1 \  \ \text{and} \ \ \EE \| \mathbb{G}^X_n \|^2_{\mathcal{G}_\delta}\ \lesssim \delta^2.
\end{align*}
\end{lemma}

The next lemma is a slight extension of Proposition 14 in \cite{azadkia2022linear} to the case where $\mathcal B_0$ is not necessarily reduced to a single point. The proof is also similar except that it uses Lemma \ref{lem: entropy} above instead of Proposition 11 in \cite{azadkia2022linear}. Hence, the proof will be omitted.
Let us define the population criterion as  
\begin{align}\label{eq: defcalD} 
\mathcal{D}(\beta) :=  \int\Big(F^Y(y) - \int F^\epsilon(y-\beta^\top  x)d F^X(x) \Big)^2 d F^Y(y)
\end{align}
where  $F^V$ is the distribution function of $V$, with $V$ denoting either $X$, or $Y$, or $\epsilon$.

\begin{lemma}\label{lem: Esup}
Suppose that $\mathbb E\| X\|^2 < \infty$,  and that $\epsilon$ has a bounded density.
Then, there exists a constant $C>0$, such that
\begin{align*}
\sqrt n\EE\sup_{\beta,\beta_0}\Big| \DD_n(\beta)-\mathcal D(\beta)-(\DD_n(\beta_0)-\mathcal D(\beta_0))\Big|  \leq C\Big(\delta+\frac\delta{\sqrt n}+\frac 1{\sqrt n}\Big),
\end{align*}
where the supremum is extended to all $\beta\in\RR^d$ and $\beta_0\in\mathcal B_0$ such that $\|\beta-\beta_0\|\leq\delta.$
\end{lemma}

The next lemma will be useful to bound the Wasserstein distance between $\mu_0$ and its estimator, it uses the theory of entropy with bracketing. For a class of functions $\mathcal F$ on $\RR^d$ and $\varepsilon>0$  we denote by $ N_{[]}(\varepsilon,\mathcal F)$  the bracketing number of $\mathcal F$ with radius $\varepsilon$, with respect to the $L^2( \PP^X)$-norm, see Definition 2.1.6 in \cite{vdvwellner23}.

\begin{lemma}\label{lem: entropy2}
For every $\beta_0\in\mathcal B_0$ and $ t\in\RR$,  consider the function defined for $x \in\RR^d$ by 
\begin{align*}
f_{\beta_0,t}({x})=\begin{cases}\one_{\beta_0^\top  x>t}&\mbox{if } t>0\\
\one_{\beta_0^\top  x\leq t}&\mbox{if }t<0.
\end{cases}
\end{align*}
 and let $\mathcal F_t$ be the class of functions ${\mathcal F_t}=\{f_{\beta_0,{t}},\beta_0\in\mathcal B_0\}$. Suppose that  $\EE\|X\|^p<\infty$ for some $p>2$, and there exists $A>0$ such that for all intervals $[a,b]$, $\mu_0([a,b])\leq A|b-a|$. 
\begin{enumerate}
\item There exists $C>0$ such that for all  $t\in\RR$ and  $\delta>0$, 
\begin{align*}
\log N_{[]}(\delta,\mathcal F_{t})\leq 
C(1+ \log(\delta^{-1})).
\end{align*}
\item For arbitrary $q\in(0,2)$, there exists $C>0$  such that  for all ${t}$ such that $|{t}|>1$, one has
\begin{align*}
\sqrt n\EE\sup_{f\in{\mathcal F_t}}\Big|\int f(d\PP_n^X-d{\mathbb P}^X)\Big|
\leq C \Big(
|t|^{-\frac p2(1-\frac q2)}+|t|^{\frac{pq}2}n^{-1/2}\Big).
\end{align*}
\item 
There exists $C>0$  such that such that for all  ${t}$ such that $|t|\leq 1$, one has 
\begin{align*}
\sqrt n\EE\sup_{f\in{\mathcal F_t}}\Big|\int f(d\PP_n^X-d{\mathbb P}^X)\Big|\leq C.
\end{align*}
\end{enumerate}
\end{lemma}


\subsection{Proof of Theorem~\ref{theo: conv2}}  
For $\beta \in \RR^d$, let us define
\begin{align*}
C_{\beta}(y) :=  \int F^\epsilon(y - \beta^\top  x)d F^X(x) \ \  \text{and} \ \ C_{n, \beta}(y) := \int  F^\epsilon(y - \beta^\top  x)  d F_n^X(x)
\end{align*}
as the convolution distribution function and its empirical estimator, where $F_n^X$ is the empirical distribution function in the sample $\mathcal X_n$ and $F^X$ is the distribution function of $X$. Recall the definition of $\mathcal{D}$ from (\ref{eq: defcalD}).   It is proved in \cite{azadkia2022linear} that 
\begin{align}\label{eq: cvcriterion}
\sup_{\beta\in\RR^d}| \DD_{n}(\beta) - \mathcal{D}(\beta)| \stackrel{\PP}{\to} 0,
\end{align}
see the display before (A.10) in that paper. Note that although Theorem 5 in  \cite{azadkia2022linear} assumes that $\mathcal B_0$ is reduced to a single point, this assumption is not used for the proof of the convergence in the last display, whence the convergence still holds under our assumptions.
With similar arguments as for the proof of (A.10) in \cite{azadkia2022linear}, where we just replace the set $\mathcal O$ by the closed ball of center $\beta_0$ and radius $r$,  one obtains that for arbitrary $\beta_0\in\mathcal B_0$ and $r>0$,
\begin{align}\label{sep}
0= \mathcal{D}(\beta_0)  <  \inf_{\beta: \, d(\beta,\mathcal B_0)>r} \mathcal{D}(\beta).
\end{align}
Note that the corresponding proof in \cite{azadkia2022linear} uses the assumption that $\mathcal B_0$ contains a unique vector $\beta_0$ to claim that if 
\begin{align*}
\int\Big(C_{\tilde{\beta}}(y) - C_{\beta_0}(y)\Big)^2 d F^Y(y) = 0,
\end{align*}
for some $\tilde\beta$, then we must have  $\tilde{\beta} = \beta_0$.   In our case, $\mathcal B_0$ may contain several vectors so we cannot use that argument. Instead, we use the fact that if the equality in the above display holds, then $\tilde{\beta}^\top X+\epsilon$ must have the same distribution as  $ \beta_0^\top X+\epsilon$, which by assumption \eqref{eq: identifiability}  implies that $\tilde{\beta}^\top X$ has the same distribution as  $ \beta_0^\top X$, which in turn means that $\tilde\beta\in\mathcal B_0$ by definition of $\mathcal B_0$, and hence $d(\tilde\beta,\mathcal B_0)=0$. The rest of the proof is similar to that in \cite{azadkia2022linear}.
Now, for arbitrary $r>0$ it follows from \eqref{sep} that we can find  $\eta_r>0$ such that
\begin{align*}
   \inf_{\beta: \, d(\beta,\mathcal B_0)>r} \mathcal{D}(\beta)-\mathcal{D}(\beta_0)  = \inf_{\beta: \, d(\beta,\mathcal B_0)>r} \mathcal{D}(\beta)  >\eta_r.
\end{align*}
Hence, it follows from \eqref{eq: cvcriterion},  \eqref{sep} and the definition of $\widehat \beta_n$,  that for arbitrary $r>0$, arbitrary $\beta_0\in\mathcal B_0$, and $\widehat{\beta}_n\in\mathcal{B}_{n}$ one has
\begin{align*}
\limsup_{n\to\infty}\PP\Big({d}(\widehat\beta_n,\mathcal B_0)>r\Big)
&\leq
\limsup_{n\to\infty}\PP\Big( \DD_n(\beta_0)  \geq   \inf_{\beta: \,d(\beta,\mathcal B_0)>r} \DD_n(\beta)\Big)\\
&\leq
\PP\Big(o_\PP(1)>\eta_r\Big),
\end{align*}
which converges to zero as $n\to\infty$. Therefore, ${d}(\widehat\beta_n,\mathcal B_0)$ converges in probability to zero as $n\to\infty$. This proves the first claim.

We turn to the proof of the second claim. 
By the triangle inequality, it follows that for arbitrary $\beta_0\in\mathcal B_0$ one has
\begin{align*}
W_1(\mu_{\widehat\beta_n}^*,\mu_0)
\leq W_1(\mu_{\widehat\beta_n}^*,\mu_{\beta_0}^*)+W_1(\mu_{\beta_0}^*,\mu_0)
\end{align*}
and therefore,
\begin{align}\label{eq: triangleW}
W_1(\mu_{\widehat\beta_n}^*,\mu_0)
\leq \inf_{\beta_0\in\mathcal B_0} W_1(\mu_{\widehat\beta_n}^*,\mu_{\beta_0}^*)+\sup_{\beta_0\in\mathcal B_0} W_1(\mu_{\beta_0}^*,\mu_0).
\end{align}
We will consider both terms on the right-hand side separately, and show that they both converge in probability to zero as $n\to\infty$. 
To deal with the first term, note that by definition of the Wasserstein{-1} distance, we have
\begin{align*}
W_1(\mu_{\widehat\beta_n}^*,\mu_{\beta_0}^*)=\inf_{Z_1,Z_2}E|Z_1-Z_2|
\end{align*}
where the infimum is taken over  all pairs $(Z_1,Z_2)$ such that $Z_1\sim\mu_{\widehat\beta_n}^*$ and $Z_2\sim\mu_{\beta_0}^*$. Taking the particular choice 
$Z_1=\widehat\beta_n^\top X^*$ and $Z_2=\beta_0^\top X^*$, we get
\begin{align}\label{eq: majWnorm} \notag
W_1(\mu_{\widehat\beta_n}^*,\mu_{\beta_0}^*)
&\leq E|(\widehat\beta_n-\beta_0)^\top X^*|\\
&\leq \|\widehat\beta_n-\beta_0\| E\|X^*\|,
\end{align}
where in the second inequality we use the Cauchy-Schwarz inequality. Taking the infimum over all possible $\beta_0$s on both sides yields
\begin{align*}
\inf_{\beta_0\in\mathcal B_0}  W_1(\mu_{\widehat\beta_n}^*,\mu_{\beta_0}^*)
\leq{d}(\widehat\beta_n,\mathcal B_0)  E\|X^*\|
\end{align*}
The second term on the right-hand side is bounded in probability since, by definition of $X^*$, its expectation is finite:
\begin{align*}
\EE E\|X^*\| =  \EE\Big(\frac 1n \sum_{i=1}^n\|X_i\|\Big)
 =  \EE \|X\|   <\infty
\end{align*}
where we recall that $E$ denotes the conditional expectation given $\{\mathcal{X}_n,\mathcal{Y}_n\}$.
Since ${d}(\widehat\beta_n,\mathcal B_0)$ converges to zero in probability by the first claim of Theorem \ref{theo: conv2}, we get
\begin{align*}
\inf_{\beta_0\in\mathcal B_0}  W_1(\mu_{\widehat\beta_n}^*,\mu_{\beta_0}^*)
= o_\PP(1).
\end{align*}
It remains to show that for the second term on the right-hand side of \eqref{eq: triangleW}, we have
\begin{align}\label{eq: second}
\sup_{\beta_0\in\mathcal B_0} W_1(\mu_{\beta_0}^*,\mu_0)
=o_\PP(1).
\end{align}
 For all $\beta_0$, $\mu_{\beta_0}^*$ is the empirical distribution in the sample $\beta_0^\top X_1,\dots,\beta_0^\top X_n$ whereas $\mu_0$ is the common distribution of the i.i.d.\ variables $\beta_0^\top X_1,\dots,\beta_0^\top X_n$. Moreover, it is assumed that $\mu_0\in \mathcal C(M,a)$ for some $M,a>0$ so, with $F_0$ the distribution function corresponding to $\mu_0$ and $Z\sim\mu_0$, we have
\begin{align*}
\int_\RR\sqrt{F_0(t)(1-F_0(t))}dt 
&\leq
2\int_0^\infty \sqrt{\PP(|Z|\geq t)}dt\\
&\leq
2\Big(1+\int_1^\infty \sqrt{\frac{\mathbb{E}|Z|^{a+2}}{t^{a+2}}}dt\Big)\\
&\leq 
2\Big(1+\int_1^\infty \frac{M^{1/2}}{t^{1+a/2}}dt\Big)\\
&\leq C
\end{align*}
for some $C>0$ that depends only on $a$ and $M$. Hence, Theorem 3.2 in \cite{bobkov2019one} implies that
\begin{align*}
\sup_{\beta_0\in\mathcal B_0}\EE W_1(\mu_{\beta_0}^*,\mu_{0})\leq Cn^{-1/2}.
\end{align*}
Let $\eta\in(0,1)$ to be fixed later, $N$ be the minimal number of balls in $\RR^d$ with radius $\eta$ that cover $\mathcal B_0$, and $\beta_1,\dots,\beta_{N}$ be the centers of such balls. 
Lemma \ref{lem: projection} ensures that $\mathcal B_0$ is a compact subset of $\RR^d$ so it follows from Lemma 2.5 in \cite{van2000empirical} that we can find $K>0$ that depends only on the radius of $\mathcal B_0$ and $d$ such that $N\leq K\eta^{-d}$, and we have
\begin{align*}
\sup_{\beta_0\in\mathcal B_0} W_1(\mu_{\beta_0}^*,\mu_0)
\leq
\max_{1\leq j\leq N}\Big(W_1(\mu_{\beta_j}^*,\mu_0)+ \sup_{\beta\in\RR^d:\ \|\beta-\beta_j\|\leq \eta}W_1(\mu_{\beta}^*,\mu_{\beta_j}^*)\Big).
\end{align*}
Note that for given $j$, $\beta_j$ need not to belong to $\mathcal B_0$. However, because the $\beta_j$s define a minimal $\eta$-covering of $\mathcal B_0$, one can find $\theta_j\in\mathcal B_0$ such that $\|\beta_j-\theta_j\|\leq \eta$ which implies that
\begin{align*}
\sup_{\beta_0\in\mathcal B_0} W_1(\mu_{\beta_0}^*,\mu_0)
&\leq
\max_{1\leq j\leq N}\Big(W_1(\mu_{\theta_j}^*,\mu_0)+ 2\sup_{\beta\in\RR^d:\ \|\beta-\beta_j\|\leq \eta}W_1(\mu_{\beta}^*,\mu_{\beta_j}^*)\Big)\\
&\leq
\sum_{j=1}^N W_1(\mu_{\theta_j}^*,\mu_0)+ 2\sup_{\beta,\beta'\in\RR^d:\ \|\beta-\beta'\|\leq \eta}W_1(\mu_{\beta}^*,\mu_{\beta'}^*).
\end{align*}
Similar to \eqref{eq: majWnorm} one obtains
\begin{align*}
\sup_{\beta,\beta'\in\RR^d:\ \|\beta-\beta'\|\leq \eta}W_1(\mu_{\beta}^*,\mu_{\beta'}^*)
\leq \eta \|X^*\|
\end{align*}
whence
\begin{align*}
\EE \sup_{\beta,\beta'\in\RR^d:\ \|\beta-\beta'\|\leq \eta}W_1(\mu_{\beta}^*,\mu_{\beta'}^*)
\leq \eta \EE \|X^*\| = \eta \EE\|X\|.
\end{align*}
Combining, we get 
\begin{align*}
\EE \sup_{\beta_0\in\mathcal B_0} W_1(\mu_{\beta_0}^*,\mu_0)
&\leq
N\sup_{\beta_0\in\mathcal B_0}\EE W_1(\mu_{\beta_0}^*,\mu_{0})+ 2\eta \EE \|X\|\\
&\leq C(\eta^{-d} n^{-1/2}+ \eta)
\end{align*}
for some constant $C>0$. We can choose $\eta$ in such a way that $\eta\ll 1$ and $\eta\gg n^{-1/(2d)}$ so that the right-hand side converges to zero as $n\to\infty$. For such a choice of $\eta$ we get 
 \eqref{eq: second}, which concludes the proof of Theorem \ref{theo: conv2}. \hfill{$\Box$}


\subsection{Proof of Theorem~\ref{theo: rateCV}}
Let us recall again the definition of $\mathcal D$ (aleady given in (\ref{eq: defcalD})): 
\begin{align*}
\mathcal{D}(\beta) :=  \int\Big(F^Y(y) - \int F^\epsilon(y-\beta^\top  x)d F^X(x) \Big)^2 d F^Y(y).
\end{align*}
By Lemma \ref{lem: projection}, we know that $\mathcal B_0$ is compact. This implies that for arbitrary $\eta>0$,  the set $\mathcal C_\eta$ of all $\beta\in\RR^d$ such that ${d}(\beta,\mathcal B_0)\leq\eta$ is a compact subset  of $\RR^d$. 
For all $\beta\in\mathcal C_\eta$,  we denote by $\beta_0$ (which depends on $\beta$) an element of $\mathcal B_0$ that satisfies $\|\beta-\beta_0\|\leq \eta$. It is clear that there exists at least one such $\beta_0$ by compactness.
In the identifiable case, it is proved in \cite{azadkia2022linear} that $\mathcal D$ is twice continuously differentiable on a small neighbourhood of (the unique) $\beta_0$, see the proof of Proposition 15 in that paper where $\mathcal D$ was denoted by $\psi$. A closer inspection of that proof shows that differentiability holds more generally on $\mathbb R^d$.
 Moreover, with similar arguments as in the proof of Proposition 15  in \cite{azadkia2022linear}, we obtain that for all $\beta_0\in\mathcal B_0$,
\begin{align*}
\frac{\partial^2\mathcal D(\beta)}{\partial\beta\partial\beta^\top }\Big|_{\beta=\beta_0}=2U(\beta_0),
\end{align*}
and that second derivative is uniformly continuous on the compact set ${\mathcal C_\eta}$. Now, $\mathcal D(\beta)\geq 0$ for all $\beta\in\RR^d$ with  exact equality if $\beta\in\mathcal B_0$.  Thus, the first derivative of $\mathcal D$ equals zero on $\mathcal B_0$. Using Taylor expansion, it follows that
\begin{align*}
\mathcal D(\beta) &= (\beta-\beta_0)^\top U(\beta_0)(\beta-\beta_0)+o(\|\beta-\beta_0\|^2)\\
&\geq\|\beta-\beta_0\|^2 \lambda(\beta_0)+o(\|\beta-\beta_0\|^2)\\
&\geq\|\beta-\beta_0\|^2 \lambda_0+o(\|\beta-\beta_0\|^2)
\end{align*}
where
\begin{align*}
\lambda_0 =\inf_{\beta_0 \in \mathcal{B}_0} \lambda(\beta_0) , 
\end{align*}
which is strictly positive by assumption.
The small $o$-term is uniform in $\beta$ since $U$ is continuous, and whence uniformly continuous on the compact set $\mathcal C_\eta$. This implies that  we can choose $\eta$ small enough so that
\begin{align*}
o(\|\beta-\beta_0\|^2)\geq -\lambda_0\|\beta-\beta_0\|^2/2
\end{align*}
for all $\beta\in\mathcal C_\eta$. Hence,
\begin{align*}
\mathcal D(\beta)\geq  \lambda_0\|\beta-\beta_0\|^2/2
\end{align*}
for all $\beta\in\mathcal C_\eta$, and
\begin{align*}
\mathcal D(\beta)\geq \lambda_0 \inf_{\beta_0\in\mathcal B_0}\|\beta-\beta_0\|^2/2.
\end{align*}
Next, we define $r_n=\sqrt n$ and
\begin{align*}
\phi_n(\delta)=\delta+\frac\delta{\sqrt n}+\frac 1{\sqrt n}
\end{align*}
for all $\delta>0$. Thanks to Lemma \ref{lem: Esup} and Theorem \ref{theo: conv2}, the conditions of Lemma \ref{lem: rateofCV} are satisfied with $\MM_n=-\DD_n$, $\MM=-\mathcal D$, $C_0=\sqrt{\lambda_0}$, $\alpha=1$ and $\mathcal B=\mathcal B_0$, and  the first claim in Theorem~\ref{theo: rateCV} follows. 

We turn to the proof of the second claim. For this task, we fix  $\varepsilon>0$ arbitrarily small and we note that from the first claim in Theorem~\ref{theo: rateCV}  we can find a constant $C_\varepsilon>0$ such that
\begin{align*}
\PP\Big( \sqrt {n\lambda_0} {d}(\widehat\beta_n,\mathcal B_0)\geq C_\varepsilon\Big)\leq\varepsilon/2.
\end{align*}
Therefore, for every $C>0$ we have
\begin{align}\label{eq: Pdist}
\notag
\PP\Big(\sqrt nW_1(\mu^*_{\widehat\beta_n},\mu_0)>C\Big)
&\leq
\PP\Big(\sqrt nW_1(\mu^*_{\widehat\beta_n},\mu_0)>C\mbox{ and }\sqrt {n\lambda_0} {d}(\widehat\beta_n,\mathcal B_0)< C_\varepsilon\Big)+\varepsilon/2\\
&\leq
\PP\Big(\sqrt n\sup_{\beta:\, \sqrt {n\lambda_0} {d}(\beta,\mathcal B_0)< C_\varepsilon }W_1(\mu^*_{\beta},\mu_0)>C\Big)+\varepsilon/2.
\end{align}
For all $\beta$ such that $\sqrt {n\lambda_0} {d}(\beta,\mathcal B_0)< C_\varepsilon$ we denote by $\beta_0$ (which depends on $\beta$) an element of $\mathcal B_0$ that satisfies $\sqrt {n\lambda_0}\|\beta-\beta_0\|\leq C_\varepsilon$. It is clear that there exists at least one such $\beta_0$. The triangle inequality yields 
\begin{align*}
W_1(\mu^*_{\beta},\mu_0)
\leq
W_1(\mu^*_{\beta},\mu^*_{\beta_0})+W_1(\mu^*_{\beta_0},\mu_{0})
\end{align*}
where we recall that $\mu_{\beta_0}^*$ is the distribution (given $X_1,\dots,X_n$)  of $\beta_0^\top X_{J_n}$ where $J_n$ is independent of the observations and follows a uniform distribution on $\{1,\dots,n\}$. Therefore,
\begin{align*}
\sup_\beta W_1(\mu^*_{\beta},\mu_0)
\leq
\sup_\beta W_1(\mu^*_{\beta},\mu^*_{\beta_0})+\sup_\beta W_1(\mu^*_{\beta_0},\mu_{0})
\end{align*}
where all suprema are taken over the set of $\beta$s such that $\sqrt {n\lambda_0} d(\beta,\mathcal B_0)< C_\varepsilon $.
By definition of the Wasserstein distance, for the first term on the right-hand side, we have
\begin{align*}
\sup_\beta W_1(\mu^*_{\beta},\mu^*_{\beta_0})
&\leq \sup_\beta  E| \beta^\top X_{J_n}-\beta_0^\top X_{J_n}|\\
&\leq
\sup_\beta \|\beta-\beta_0\|E\|X_{J_n}\|
\end{align*}
where $E$ is the expectation with respect to $J_n$ and 
using the Cauchy-Schwarz inequality for the second inequality. Letting 
\begin{align*}
m_n= \frac1n\sum_{i=1}^n \|X_i\|,
\end{align*}
we have $E\|X_{J_n}\|=m_n$ and therefore,
\begin{align*}
\sqrt n \sup_\beta W_1(\mu^*_{\beta},\mu^*_{\beta_0})
\leq
C_\varepsilon m_n\lambda_0^{-1/2}.
\end{align*}
Note that $\EE m_n=\EE\|X\|<\infty$ by assumption, thus, $m_n$ is bounded in probability and we get
\begin{align}\label{eq: sup+}\notag
\sqrt n \sup_\beta W_1(\mu^*_{\beta},\mu_0)
&\leq \lambda_0^{-1/2} O_\PP(1)+\sqrt n\sup_\beta W_1(\mu^*_{\beta_0},\mu_{0})   \\
& \leq   \lambda_0^{-1/2}{O_\PP(1)+\sqrt n\sup_{\beta_0 \in \mathcal{B}_0} W_1(\mu^*_{\beta_0},\mu_{0}) }.
\end{align}
 Now, for the second term on the right-hand side, we have for arbitrary $\beta_0 \in \mathcal{B}_0$
\begin{align*}
W_1(\mu^*_{\beta_0},\mu_{0})=\int_\RR|F_{\beta_0}^*(t)-F_0(t)|dt
\end{align*}
where $F_{\beta_0}^*$ and $F_0$ are the distribution functions of $\mu^*_{\beta_0}$ and $\mu_{0}$ respectively, see \cite[Theorem 2.9]{bobkov2019one}. We have \begin{align*}
F_0(t)=P(\beta_0^\top X\leq t)=1-P(\beta_0^\top X> t)
\end{align*} 
with $X$ independent of the observations and
\begin{align*}
F_{\beta_0}^*(t)=\frac 1n\sum_{i=1}^n \one_{\beta_0^\top X_i\leq t}=1-\frac 1n\sum_{i=1}^n \one_{\beta_0^\top X_i> t}.
\end{align*}
Let $\PP_n^X$ be the empirical distribution of the sample $X_1,\dots,X_n$ and ${\PP}^X$ the distribution of $X$, for all $\beta, {x}$, $f_{\beta, {t}}({x})=\one_{\beta^\top  {x}> {t}}$ if ${t}>0$ and  $f_{\beta,{t}}({x})=\one_{\beta^\top  {x}\leq {t}}$ if ${t}\leq 0$, we have
\begin{align*}
| F_{\beta_0}^*(t)-F_0(t)|=\Big|\int f_{\beta_0, {t}}d\PP_n^X-\int f_{\beta_0,  t }d \PP^X\Big|.
\end{align*}
Combining yields 
\begin{align*}
\sqrt n\sup_{\beta_0\in\mathcal B_0} W_1(\mu^*_{\beta_0},\mu_{0})
\leq\sqrt n\int_\RR\sup_{\beta_0\in\mathcal B_0}\Big|\int f_{\beta_0,  t }(d\PP_n^X-d \PP^X)\Big|d  t .
\end{align*}
Hence, taking expectations
\begin{align*}
\sqrt n\EE \sup_{\beta_0\in\mathcal B_0} W_1(\mu^*_{\beta_0},\mu_{0})
\leq E_1+E_2+E_3
\end{align*}
where for $i\in\{1,2,3\}$ we define
\begin{align*}
E_i:=
\sqrt n\int_{D_i}\EE\sup_{\beta_0\in\mathcal B_0}\Big|\int f_{\beta_0,t}(d\PP_n^X-d \PP^X)\Big|dt
\end{align*}
with the domains of integration being defined by
\begin{align*}
D_1&:=[-1,1]\\
D_2&:=\{  t :\ |  t |\in(1,a_n]\}\\
D_3&:=\{  t :\ |  t |>a_n\}
\end{align*}
where $a_n$ is a sequence that tends to infinity as $n\to\infty$ and that will be choosen later.
We deal separately with the three terms $E_i$. For $E_1$ we have
\begin{align*}
E_1\leq 2\sqrt n\sup_{  t \in[-1,1]}\EE \sup_{\beta_0\in\mathcal B_0}\Big|\int f_{\beta_0,  t }(d\PP_n^X-d \PP^X)\Big|,
\end{align*}
which is finite thanks to the third claim in Lemma \ref{lem: entropy2}. The second  claim of the same lemma proves that for $E_2$ we have
\begin{align*}
E_2&\leq C \int_{|  t |\in(1,a_n]}     \Big( |t|^{-\frac p2(1-\frac q2)}+ |  t |^{\frac{pq}2}n^{-1/2}
\Big)dt.
\end{align*}
Thus, choosing $q\in(0,2)$ small enough so that $\frac p2(1-\frac q2)>1$ and taking $a_n=n^{1/(2+pq)}$ we get that
\begin{align*}
E_2&\leq 2C \Big(\frac{a_n^{1-\frac p2(1-\frac q2)}}{1-\frac p2(1-\frac q2)}+\frac{a_n^{1+\frac{pq}2}}{1+\frac{pq}2}\,\, n^{-1/2}\Big)
\end{align*}
is also finite.

Finally, we deal with  $E_3$ with the above choice for $a_n$. We recall that $\mathcal B_0$ is compact (see Lemma  \ref{lem: projection}) and therefore, we can find $B>0$ such that $\|\beta_0\|< B$ for all $\beta_0\in\mathcal B_0$ and also, $|\beta_0^\top t|<B\|t\|$ for all $t$, using the Cauchy-Schwarz inequality. This implies that
\begin{align*}
E_3&\leq 
\sqrt n\int_{|t|>a_n}\EE\sup_{\beta_0\in\mathcal B_0}\Big(\int f_{\beta_0,t}d\PP_n^X+\int f_{\beta_0,t}d \PP^X\Big)dt\\
&\leq 
\sqrt n\int_{|t|>a_n}\EE\Big(\int \one_{B\|x\|> |t|}d\PP_n^X+\int \one_{B\|x\|> |t|}d \PP^X\Big)dt\\
&= 2\sqrt n \int_{|t|>a_n}\PP\Big(B\|X\|>|t|\Big) dt.
\end{align*}
For $p>2$ such that $\EE\|X\|^p<\infty$, by the Markov inequality it follows that there exists a constant $C>0$ such that
\begin{align*}
E_3&\leq   
C\sqrt n \int_{|t|>a_n} |t|^{-p}dt=2C\sqrt n \,\, \frac{a_n^{1-p}}{p-1}.
\end{align*}
Since $a_n=n^{1/(2+pq)}$ for some $p>2$ and a sufficiently small $q>0$, we can choose $q>0$ small enough so that $\sqrt na_n^{1-p}$ tends to zero as $n\to\infty$. Hence,  $E_i$ is finite for all $i\in\{1,2,3\}$ and therefore, as a consequence of the Markov inequality we obtain
\begin{align*}
\sqrt n\sup_{\beta_0\in\mathcal B_0} W_1(\mu^*_{\beta_0},\mu_{0})=O_\PP(1).
\end{align*}
Hence, from \eqref{eq: sup+} it follows that 
\begin{align*}
\sqrt n \sup_{\beta} W_1(\mu^*_{\beta},\mu_0)
=O_\PP(1).
\end{align*}
Thus, we can find $C>0$ large enough so that 
\begin{align*}
\PP\Big(\sqrt n\sup_\beta W_1(\mu^*_{\beta},\mu_0)>C\Big)\leq\varepsilon/2,
\end{align*}
and from \eqref{eq: Pdist} it follows that for arbitrary $\varepsilon>0$ we can find $C>0$ such that
\begin{align*}
\PP\Big(\sqrt nW_1(\mu^*_{\widehat\beta_n},\mu_0)>C\Big)
\leq\varepsilon.
\end{align*}
This completes the proof of Theorem \ref{theo: rateCV}. \hfill{$\Box$}


\subsection{Proof of Theorem \ref{theo: kde}}
    First note that 
\begin{align*}
\|\widehat f_Z-f_Z\|_{\infty} &\leq \|\mu^*_{\widehat\beta_n}*K_h - \mu_0* K_h\|_{\infty} + \|\mu_0*K_h -f_Z\|_{\infty}\\
&=:A+B.
\end{align*}
    The term $A$ can be handled via Wasserstein duality. Indeed, setting $\phi_z(y) = K_h(z-y)$, one has
    \begin{align*}
    (\mu^*_{\widehat\beta_n}*K_h - \mu_0* K_h)(z) &= \int \phi_z d(\mu^*_{\widehat\beta_n}-\mu_0) \\
    &\leq \operatorname{Lip} (\phi_z)W_1(\mu^*_{\widehat\beta_n},\mu_0) \\
    &\leq \frac{L}{h^2}W_1(\mu^*_{\widehat\beta_n},\mu_0)\\
    &= O_{\mathbb{P}}(h^{-2}n^{-1/2}),
        \end{align*}
    where we used that $W_1(\mu^*_{\widehat\beta_n},\mu_0)= O_{\mathbb{P}}(n^{-1/2})$ by Theorem \ref{theo: rateCV}, and that the kernel function is $L$-Lipschitz implying
    $$
    |\phi_z(y_1) - \phi_z(y_2)|= \frac{1}{h} \Big| K\Big(\frac{z-y_1}{h}\Big) - K\Big(\frac{z-y_2}{h}\Big) \Big| \leq \frac{L}{h^2} |y_1-y_2|
    $$
    for all $y_1,y_2\in \RR$, thus, yielding $\operatorname{Lip} (\phi_z)=  \sup_{y_1 \neq y_2} \frac{|\phi_z(y_1)-\phi_z(y_2)|}{|y_1-y_2|}\leq \frac{L}{h^2}$.
    The term $B$ is a bias term and to bound it we rewrite 
    \begin{align*}
    (\mu_0*K_h -f_Z)(z) &= \int K_h(z-y) f_Z(y) dy -f_Z(z) \\
    &= \int \frac{1}{h} K\left(\frac{z-y}{h}\right) (f_Z(y) - f_Z(z)) dy \\
    &=\int K(u) (f_Z(z-hu)-f_Z(z)) du
    \end{align*}
    where the third step is a change of variable $u=(z-y)/h$. 
    Now, by the assumptions made on $f_Z$, Taylor's Theorem implies that for all $z,u$ there exists some $\theta \in (0,1)$ such that
    $$f_Z(z - h u) - f_Z(z) =  - h u f_Z'(z) + \frac{(h u)^2}{2} f_Z''(z - \theta h u) .$$
    Since $\int uK(u)du=0$, we obtain that for all $z\in \mathbb{R}$,
    \begin{align*}
        |(\mu_0*K_h -f_Z)(z)| &=\left| \int K(u) \frac{(h u)^2}{2} f_Z''(z - \theta h u)du\right| \\
        &\leq  C h^2 \int u^2 K(u) du = O(h^2),
     \end{align*} 
    for some $C>0$. 
    Combining the bounds and since $h=h_n\asymp n^{-1/8}$ we conclude
    \begin{align}\label{rateforfZ}
    \|\widehat f_Z-f_Z\|_{\infty} = O_\mathbb{P}(n^{-1/4}).
    \end{align}
    Now, we turn to the second claim. First, let us consider the decomposition
    \begin{align}\label{identity}
	\widehat{f}_{Z|Y}(z|y_0) -{f}_{Z|Y}(z|y_0) &= \frac{f_\epsilon(y_0 - z) \widehat{f}_Z(z)}{\widehat{f}_Y(y_0)}-  \frac{f_\epsilon(y_0 - z) {f}_Z(z)}{{f}_Y(y_0)} \notag \\
    &=f_\epsilon(y_0 - z) \Big( \frac{\widehat{f}_Z(z)}{\widehat{f}_Y(y_0)} - \frac{{f}_Z(z)}{\widehat{f}_Y(y_0)} +\frac{{f}_Z(z)}{\widehat{f}_Y(y_0)} - \frac{{f}_Z(z)}{{f}_Y(y_0)} \Big) \notag \\
    &=f_\epsilon(y_0 - z) \Big( \frac{\widehat{f}_Z(z)-{f}_Z(z)}{\widehat{f}_Y(y_0)}  +{f}_Z(z)(\frac{1}{\widehat{f}_Y(y_0)} - \frac{1}{{f}_Y(y_0)} )\Big) \notag\\
    &=f_\epsilon(y_0 - z) \Big( \frac{\widehat{f}_Z(z)-{f}_Z(z)}{\widehat{f}_Y(y_0)}  - \frac{f_Z(z)}{f_Y(y_0)\widehat f_Y(y_0) } (\widehat {f}_{Y}(y_0)-f_Y(y_0))\Big).
     \end{align}
     Now,  recall that 
     \begin{align*}
 \widehat f_Y(y)  =  \frac{1}{n} \sum_{i=1}^n f_\epsilon(y- \widehat \beta^\top_n X_i)  =  \int f_\epsilon(y- z) d \mu^*_{\widehat \beta_n}(z)
     \end{align*}
 for $y \in \mathbb R$. Then, 
    \begin{align*}
 \vert \widehat f_Y(y)  - f_Y(y)  \vert  & =   \left \vert \int f_\epsilon(y- z) d (\mu^*_{\widehat \beta_n}(z) - \mu_0(z))  \right \vert \notag  \\
 & \le   M W_1(\mu^*_{\widehat \beta_n}, \mu_0)  = O_{\mathbb P}(n^{-1/2}),
     \end{align*}
and hence
\begin{align}\label{rateforfY}
    \Vert \widehat f_Y  - f_Y  \Vert_\infty = O_{\mathbb P}(n^{-1/2}).
\end{align}
Thus, with probability tending to 1, $\widehat f_Y(y_0) >  f_Y(y_0)/2$, and it follows from (\ref{rateforfZ}) and (\ref{identity}) that 
\begin{align*}
\|\widehat{f}_{Z|Y}(\cdot|y_0)-f_{Z|Y}(\cdot|y_0)\|_{\infty}  & \le   \Vert f_\epsilon \Vert_\infty   \Big( \frac{2 \Vert \widehat f_Z - f_Z \Vert_\infty}{f_Y(y_0)}  +  \frac{2 \Vert f_Z \Vert_\infty }{f_Y(y_0)^2}  \Vert \widehat f_Y - f_Y \Vert_\infty    \Big) \\
& = O_{\mathbb{P}}(n^{-1/4})+O_{\mathbb{P}}(n^{-1/2})  =  O_{\mathbb{P}}(n^{-1/4}),
\end{align*}
which concludes the proof.
\hfill{$\Box$}


\subsection{Proof of Proposition \ref{prop: mean-mode-quantile}}
    For the first part, it follows from the convergence rate in (\ref{rateforfY}) that  with probability tending to 1 
\begin{align*}
\widehat f_Y(y_0)  \ge f_Y(y_0)/2.
\end{align*}
Now, using the identity in (\ref{identity}) we can write that
\begin{align*}
 \big \vert  \widehat{E}(Z|Y=y_0) -  \mathbb E(Z | Y= y_0)  \big \vert 
 & =     \bigg \vert \frac{1}{\widehat f_Y(y_0)} \int z f_\epsilon(y_0 - z) (\widehat f_Z(z) - f_Z(z)) dz  \\
 & \quad -    \frac{\widehat f_Y(y_0) - f_Y(y_0) }{\widehat f_Y(y_0) f_Y(y_0)}\int  z f_\epsilon(y_0-z) f_Z(z) dz \bigg \vert \\
 & \le      \frac{2 \Vert \widehat f_Z - f_Z \Vert_\infty}{f_Y(y_0)} \int \vert z \vert  f_\epsilon(y_0 - z)  dz  \\
  & \quad +  \frac{ 2  \Vert f_Z \Vert_\infty \Vert \widehat {f}_{Y}-f_Y \Vert_\infty}{f_Y(y_0)^2}  \int \vert z \vert  f_\epsilon(y_0 - z)  dz  \\
&\le    2\Big(\frac{\Vert \widehat f_Z - f_Z \Vert_\infty} {f_Y(y_0)} + \  \frac{\Vert f_Z \Vert_\infty  \Vert \widehat f_Y - f_Y \Vert_\infty }{f_Y(y_0)^2}  \Big) \  (\vert y_0 \vert  +  \mathbb E[\vert \epsilon\vert ])
  \end{align*}
with probability tending to 1. Note that above we used  the fact that $\int \vert z \vert f_\epsilon(y_0 - z) dz =  \int \vert z - y_0 + y_0 \vert f_\epsilon(y_0 - z) dz  \le \mathbb E[ \vert \epsilon \vert]  + \vert y_0 \vert $.
Then, it follows from Theorem \ref{theo: kde} that 
\begin{align*}
\big \vert  \widehat{E}(Z|Y=y_0) -  \mathbb E(Z | Y= y_0)  \big \vert  =  O_{\mathbb P}(n^{-1/4}).
\end{align*}
     
    For the second claim, let us write $\widehat z^* : =  \widehat M(Z| Y=y_0)$. We first prove that $ \widehat z^* \overset{\PP}{\to} z^*$. For some constant $K > 0$, consider the event $\{ | \widehat z^*  - z^* | >  K \}$.  When this event occurs, we have by assumption that there exists some $\delta > 0$ (depending on $K$) such that 
    $$
    f_{Z|Y}(z^* | y_0)  - f_{Z|Y}(\widehat z^* | y_0) \ge \delta.
    $$
    On the other hand, it holds that 
    \begin{align}\label{eq: IneqSup}
     \vert \sup_{z}  \widehat f_{Z|Y}(z | y_0) -  \sup_{z} f_{Z|Y}(z | y_0) \vert \le \Vert \widehat f_{Z|Y}(\cdot | y_0)  -  f_{Z|Y}(\cdot | y_0) \Vert_\infty   
    \end{align}
In fact, we have that for all $z$,
\begin{align*}
\widehat f_{Z| Y}(z| y_0)  &\le    f_{Z| Y}(z| y_0) + \Vert \widehat f_{Z|Y}(\cdot | y_0)  -  f_{Z|Y}(\cdot | y_0) \Vert_\infty \\
& \le  \sup_{z} f_{Z|Y}(z| y_0) + \Vert \widehat f_{Z|Y}(\cdot | y_0)  -  f_{Z|Y}(\cdot | y_0) \Vert_\infty,
\end{align*}
implying that $ \sup_{z} \widehat f_{Z| Y}(z| y_0) \le \sup_{z} f_{Z| Y}(z| y_0) + \Vert \widehat f_{Z|Y}(\cdot | y_0)  -  f_{Z|Y}(\cdot | y_0) \Vert_\infty$. The same argument can be used to obtain the inequality in the other direction.  Then, Theorem \ref{theo: kde} implies that  
    \begin{align*}
     \widehat f_{Z| Y}(\widehat z^*|y_0)  = \sup_{z}  \widehat f_{Z| Y}(z|y_0)  \overset{\PP}{\to} \sup_{z} f_{Z| Y}(z|y_0) = f_{Z| Y}(z^*| y_0),
    \end{align*}
    and hence with probability tending to 1,
    \begin{align*}
    \widehat f_{Z|Y}(\widehat z^*|y_0)  +  \frac{\delta}{2} \ge  f_{Z|Y}(z^*|y_0) \ge  f_{Z|Y}(\widehat z^*|y_0)  +  \delta,
    \end{align*}
    from which we conclude that 
    $$
    \widehat f_{Z|Y}(\widehat z^*|y_0)  \ge f_{Z|Y}(\widehat z^*|y_0) + \frac{\delta}{2}
    $$
    with probability 1. Thus,
    \begin{align*}
     \mathbb P( | \widehat z^*  - z^* | >  K )  & \le \mathbb P \Big( \widehat f_{Z|Y}(\widehat z^*|y_0) -   f_{Z|Y}(\widehat z^*|y_0) \ge \frac{\delta}{2}\Big)  + o_{\mathbb P}(1) \\
     & \le   \mathbb P \Big( \Vert \widehat f_{Z|Y}(\cdot| y_0)  -  f_{Z|Y}(\cdot| y_0)  \Vert_\infty \ge \frac{\delta}{2}\Big)   + o_{\mathbb P}(1)\to 0
    \end{align*}
as $n \to \infty$. Since $K$ was arbitrary, this proves that $\widehat z^*$ is consistent. This means that for $\eta > 0$ small, 
\begin{align*}
\mathbb P(\vert \widehat z^* - z^*  \vert  \le \eta)  \to 1
\end{align*}
as $n \to \infty$. When $\eta$ is small enough, and the event  $\{\vert \widehat z^* - z^*  \vert  \le \eta \}$ occurs, then $\widehat z^*$ belongs to the small neighbourhood mentioned in the assumptions of the proposition. Fix $C > 0$ and note that
\begin{align*}
 \mathbb P(\vert \widehat z^* - z^* \vert > C n^{-1/(4\alpha)} )  \le \mathbb P (  f_{Z|Y}(\widehat z^*|y_0) -  f_{Z|Y}(z^*|y_0) \ge \lambda C^\alpha n^{-1/4} ) + o_{\mathbb P}(1). 
\end{align*}
Moreover, it holds that
\begin{align*}
 f_{Z|Y}(\widehat z^*|y_0) -  f_{Z|Y}(z^*|y_0) & =   f_{Z|Y}(\widehat z^*|y_0)  - \widehat f_{Z| Y}(\widehat z^*|y_0) +  \widehat f_{Z| Y}(\widehat z^*|y_0)- f_{Z|Y}(z^*|y_0) \\
& \le   \Vert \widehat f_{Z|Y}(\cdot | y_0)  -  f_{Z|Y}(\cdot | y_0) \Vert_\infty +  \vert \sup_{z}  \widehat f_{Z|Y}(\cdot | y_0) -  \sup_{z} f_{Z|Y}(\cdot | y_0) \vert  \\
& \le   2 \Vert \widehat f_{Z|Y}(\cdot | y_0)  -  f_{Z|Y}(\cdot | y_0) \Vert_\infty  = O_{\mathbb P}(n^{-1/4}),
\end{align*}
where we used \eqref{eq: IneqSup} in the last step. Therefore, it follows
\begin{align*}
 \mathbb P(\vert \widehat z^* - z^* \vert > C n^{-1/(4\alpha)} )  \le \mathbb P ( \Vert \widehat f_{Z|Y}(\cdot | y_0)  -  f_{Z|Y}(\cdot | y_0) \Vert_\infty  \ge \lambda C^\alpha n^{-1/4} /2 )  + o_{\mathbb P}(1) \to 0  
 \end{align*}
as $C \to \infty$ and $n \to \infty$. This proves that $\vert \widehat z^* - z^* \vert = O_{\mathbb P}(n^{-1/(4\alpha)})$.

For the third claim, let us write again $\widehat F_{Z|Y}(z|y_0):= \int_{-\infty}^z \widehat f_{Z|Y}(x|y_0)dx $.
Note that 
\begin{align*}
\sup_{z\in \mathbb R}| \widehat F_{Z|Y}(z|y_0)-F_{Z|Y}(z|y_0)|&\leq \sup_{z\in \mathbb R }\int_{-\infty}^z| \widehat f_{Z|Y}(x|y_0)-f_{Z|Y}(x|y_0)|dx\\
&\leq \int_{\RR} | \widehat f_{Z|Y}(z|y_0)-f_{Z|Y}(z|y_0)|dz \\
& =  \int_{\vert z \vert \le M} | \widehat f_{Z|Y}(z|y_0)-f_{Z|Y}(z|y_0)|dz  \\
&\quad + \int_{\vert z \vert > M} | \widehat f_{Z|Y}(z|y_0)-f_{Z|Y}(z|y_0)|dz \\
& = I_1  + I_2
\end{align*}
For the first term $I_1$, we have that
\begin{align*}
I_1 
&\leq 2M \|\widehat f_{Z|Y}(\cdot|y_0)-f_{Z|Y}.(\cdot|y_0)\|_{\infty}. 
\end{align*}
For the second term, we have for $M > c_1$ that
\begin{align*}
I_2 &= 
\int_{|z|> M}| \widehat f_{Z|Y}(z|y_0)-f_{Z|Y}(z|y_0)|dz \\
& \le \int_{|z|> M}  \widehat f_{Z|Y}(z|y_0) dz + \int_{|z|> M} f_{Z|Y}(z|y_0) dz  \\ 
& = 1 - \int_{|z| \le M}  \widehat f_{Z|Y}(z|y_0) dz  +  \int_{|z|> M} f_{Z|Y}(z|y_0) dz \\
& =  1 - \int_{|z| \le M}  \left( \widehat f_{Z|Y}(z|y_0) - f_{Z|Y}(z|y_0) dz \right) dz - \int_{|z| \le M}  f_{Z|Y}(z|y_0) dz \\
& \ \ + \int_{|z|> M} f_{Z|Y}(z|y_0) dz \\
& \le \int_{|z| \le M}   \vert \widehat f_{Z|Y}(z|y_0) - f_{Z|Y}(z|y_0)  \vert dz   + 2 \int_{|z|> M} f_{Z|Y}(z|y_0) dz \\
& \le 2M  \|\widehat f_{Z|Y}(\cdot|y_0)-f_{Z|Y}(\cdot|y_0)\|_{\infty} + 2 \int_{|z|> M} f_{Z|Y}(z|y_0) dz \\
& \le  2M  \|\widehat f_{Z|Y}(\cdot|y_0)-f_{Z|Y}(\cdot|y_0)\|_{\infty}  + \frac{4c_2}{c_3} e^{-c_3M}
\end{align*}
where the last step follows from our assumption about the tail decay of $f_{Z|Y}(\cdot|y_0)$.  Therefore, we have that 
\begin{align*}
\sup_{z\in \mathbb R}| \widehat F_{Z|Y}(z|y_0)-F_{Z|Y}(z|y_0)| & \le 4M \|\widehat f_{Z|Y}(\cdot|y_0)-f_{Z|Y}.(\cdot|y_0)\|_{\infty}  + \frac{4c_2}{c_3} e^{-c_3M}. 
\end{align*}
Now taking $M=M_n= \frac{1}{c_3}\log (\frac{4c_2}{c_3n^{-1/4}})=\frac{1}{4c_3}\log n + \frac{1}{c_3}\log ( \frac{4c_2}{c_3})$ yields $\frac{4c_2}{c_3} e^{-c_3M} = n^{-1/4}$, and using the result of Theorem \ref{theo: kde} it follows that
\begin{align}
 \label{eq: CDF-conv}
    \sup_{z\in \mathbb R }| \widehat F_{Z|Y}(z|y_0)-F_{Z|Y}(z|y_0)| = O_{\mathbb P}(n^{-1/4} \log n)  +  O_{\mathbb P}( n^{-1/4})= O_\PP(n^{-1/4} \log n).
\end{align}
Now, consider a neighbourhood $N_\delta:= [Q_{Z|Y}(p|y_0)-\delta , Q_{Z|Y}(p|y_0)+\delta]$ for some $\delta>0$.
By assumption, for any $z_1,z_2\in N_\delta$, we have
      $$|F_{Z|Y}(z_2|y_0)-F_{Z|Y}(z_1|y_0)|=\Big|\int_{z_1}^{z_2}f_{Z|Y}(x|y_0)dx\Big|\geq |z_2-z_1| \inf_{z\in N_\delta}f_{Z|Y}(z|y_0)>0.$$
Thus, since by \eqref{eq: CDF-conv} and local strict positivity of $f_{Z|Y}(\cdot|y_0)$ on $N_\delta$ we have $\widehat Q_{Z|Y}(p|y_0) \in N_\delta$ with probability tending to 1, we conclude that on this event
\begin{align*}
|\widehat Q_{Z|Y}(p|y_0)-Q_{Z|Y}(p|y_0)| &\leq \frac{|F_{Z|Y}(\widehat Q_{Z|Y}(p|y_0)|y_0)-F_{Z|Y}( Q_{Z|Y}(p|y_0)|y_0)|}{\inf_{z\in N_\delta}f_{Z|Y}(z|y_0)}\\
&\leq \frac{|F_{Z|Y}(\widehat Q_{Z|Y}(p|y_0)|y_0)-\widehat F_{Z|Y}( \widehat Q_{Z|Y}(p|y_0)|y_0)|}{\inf_{z\in N_\delta}f_{Z|Y}(z|y_0)}\\
&\qquad + \frac{|\widehat F_{Z|Y}(\widehat Q_{Z|Y}(p|y_0)|y_0)-F_{Z|Y}( Q_{Z|Y}(p|y_0)|y_0)|}{\inf_{z\in N_\delta}f_{Z|Y}(z|y_0)}\\
&\leq \frac{\sup_{z\in \mathbb{R} }|F_{Z|Y}(z|y_0)-\widehat F_{Z|Y}(  z|y_0)|}{\inf_{z\in N_\delta}f_{Z|Y}(z|y_0)}\\ &=O_{\mathbb{P}}(n^{-1/4}\log n),
\end{align*}
which concludes the proof.
\hfill{$\Box$}


\subsection{Proof of the lemmas}
We point out that our results generalize some results of \cite{azadkia2022linear} to the case of possibly non-identifiable $\beta_0$. Moreover, in the introduction of \cite{azadkia2022linear}, it is assumed that the two samples $\mathcal X_n$ and $\mathcal Y_n$ are independent, an assumption that we do not make in the present paper. However, the proofs in \cite{azadkia2022linear}, which we generalize here, do not use this assumption.

\subsubsection{{Proof of Lemma \ref{lem: identifiability}}}
For the equivalence in \eqref{eq: identifiability} to hold, it is sufficient to have the equivalence
\begin{align*}
Z +\epsilon\stackrel{d}{=} Z_0+\epsilon \ \ \Longleftrightarrow \ \ Z  \sim\mu_0 
\end{align*}
for any real random variable $Z$ which is independent of $\epsilon$. Let
\begin{align*}
\kappa_Z(t)=\EE\exp(itZ),  \ t \in \RR
\end{align*}
denote the characteristic function of a real random variable $Z$.  We recall that two given random variables have the same distribution if and only if their characteristic functions are equal. Moreover, if  $Z$ and $\epsilon$ are independent, then
\begin{align*}
\kappa_{Z+\epsilon}=\kappa_{Z}\kappa_\epsilon.
\end{align*}
Hence, a sufficient condition for identifiability of $\mu_0$ is that we have the equivalence
\begin{align*}
\kappa_{Z}\kappa_\epsilon=\kappa_{Z_0}\kappa_\epsilon \ \ \Longleftrightarrow  \ \
\kappa_{Z}=\kappa_{Z_0}
\end{align*}
for all random variables $Z$ that are independent of $\epsilon$.  
It is proved in \cite{meister2009} (see the comments below display (2.23) on page 25) that a sufficient condition for the above equivalence to hold is that the set of zeros of $\kappa_\epsilon$ does not contain any open, nonempty interval as a subset. \hfill{$\Box$}

\subsubsection{Proof of Lemma \ref{lem: rateofCV}}
Assume for simplicity that $\widehat\theta_n$ truly maximizes the map $\theta\mapsto\MM_n(\theta)$. For each $n$, $\Theta\setminus\mathcal{B}$ can be partitioned into shells 
\begin{align*}
S_{j,n}=\{\theta:\ 2^{j-1}<C_0 r_nd(\theta,\mathcal B)\leq 2^j\}
\end{align*} 
with $j\in \mathbb{N}$. Similarly to the proof of Theorem 3.2.5 in \cite{vdvwellner23} for every $\eta>0$ we get
\begin{align*}
\mathbb{P}\Big(C_0 r_nd(\widehat\theta_n,\mathcal B)>2^M\Big)
&\leq \sum_{j> M;\ 2^j\leq C_0\eta r_n}\mathbb{P}\Big(\sup_{(\theta,\theta_0)	\in S_{j,n}\times \mathcal{B}}(\MM_n(\theta)-\MM_n(\theta_0))\geq 0\Big)\\
& \quad + \mathbb{P}\Big(2d(\widehat\theta_n,\mathcal B))\geq \eta\Big).
\end{align*}
Since  $d(\widehat\theta_n,\mathcal B)$ converges in outer probability to zero, the second probability on the right converges to zero as $n\to\infty$ for every $\eta>0$. Thus, choose $\eta>0$ small enough that the condition in \eqref{eq: rateofCV} holds for every $\delta\leq\eta$. Then, for every $\theta_0\in \mathcal{B}$, for every $j$ involved in the sum and for every $\theta\in S_{j,n}$ we have
\begin{align*}
\MM(\theta)-\MM(\theta_0)\leq -C_0^2d^2(\theta,\mathcal B)< \frac{-2^{2j-2}}{r_n^2}.
\end{align*}
The rest of the proof is the same as that of \cite{vdvwellner23} and is omitted. \hfill{$\Box$}


\subsubsection{Proof of Lemma \ref{lem: projection}}\label{sec: projection}
Let $\Sigma$ denote the variance-covariance matrix of $X$ and $\beta_0\in\mathcal B_0$ be arbitrary. Then, $\mathrm{Var}(\beta_0^\top X)=\beta_0^\top \Sigma\beta_0=: c$ with $c\in(0,+\infty)$, and because $\beta^\top X$ has the same distribution as $\beta_0^\top X$ for all $\beta\in\mathcal B_0$, we have $\|\beta\|_\Sigma=c$ for all $\beta\in\mathcal B_0$, where $\|\cdot\|_\Sigma$ is the norm defined by $\|\beta\|_\Sigma=\beta^\top \Sigma\beta$. Hence, $\mathcal B_0
$ is a bounded set in $\RR^d$. 

Now, let $(\beta_n)$ be a sequence in $\mathcal B_0$ that converges to some $\beta\in\RR^d$ as $n\to\infty$. For every $n$, $\beta_n^\top X$ has the same distribution as $\beta_0^\top X$ where $\beta_0\in\mathcal B_0$ is arbitrary, hence $\beta_n^\top X$ has the same characteristic function as $\beta_0^\top X$ for every $n$. Letting $n\to\infty$ shows that $\beta^\top X$ also has the same characteristic function as $\beta_0^\top X$, whence $\beta\in\mathcal B_0$. This means that $\mathcal B_0$ is a closed subset of $\RR^d$, whence it is compact. The lemma follows. \hfill{$\Box$}


\subsubsection{Proof of Lemma \ref{lem: entropy}} 
With similar arguments as in  the proof of Proposition 11  in \cite{azadkia2022linear}, we see that  $\mathcal{G}$  is a subset of the convex hull of the class
\begin{align*}
\mathcal H:=\{x\mapsto\phi_1(\beta^\top x+t-y)\times\phi_2(\beta_0^\top x+t-y),\ \beta,\beta_0\in\RR^d,\ t,y\in\RR\},
\end{align*}
and that $\mathcal H$ is a VC-class. Here, $\phi_1(u)=\one\{u\geq 0\}$ and $\phi_2(u)=\one\{u\leq 0\}, $ $u\in\RR$.
Moreover, denoting by $M$ a positive number such that the density $f_\epsilon$ of $\epsilon$ satisfies $f_\epsilon(t)\leq M$ for all $t\in\RR$, for all $\beta,\beta'$ such that $\|\beta-\beta'\|\leq \delta$ and all $x,y$ we have
\begin{align*}
|g_{\beta,\beta', y}(x)| \leq M|(\beta-\beta')^\top  x| 
\leq M \delta \|x\|
\end{align*}
using the Cauchy-Schwarz inequality. This means that the function $G(x):= M\delta\|x\|$ is an envelope of $\mathcal G$. Since we have the same envelope as for their class $\mathcal G$, the second claim of the lemma follows from similar arguments as in the proof of Proposition 11 in \cite{azadkia2022linear}, that uses finiteness of $J(1,\mathcal G)$ and applies Theorem 2.14.1 of \cite{vdvwellner23}.

We turn to the proof of the first claim. It is proved in Section \ref{sec: projection} that $\mathcal B_0$ is a compact set in $\RR^d$. Hence, it is contained in the ball with center 0 and radius $C$, for some $C>0$.  For every $\beta_0\in\mathcal B_0$ we can write
\begin{align*}
m_{\beta_0,y}=m_{0,y}+(m_{\beta_0,y}-m_{0,y})
\end{align*}
so with
\begin{align*}
\mathcal{M}_0:= \{ m_{0, y}, y \in \RR  \},
\end{align*}
we have
\begin{align*}
\mathcal{M}\subset\mathcal{M}_0+\mathcal G(C)
\end{align*}
where $\mathcal G(C)$ is defined in the same way as $\mathcal G$ with $\delta$ replaced by $C$. Hence, we have
\begin{align*}
\EE \| \mathbb{G}^X_n \|^2_{\mathcal{M}} \leq\EE \| {\GG}^X_n \|^2_{\mathcal M_0} +\EE \| \mathbb{G}^X_n \|^2_{\mathcal G(C)} .
\end{align*}
It follows from Proposition 11  in \cite{azadkia2022linear} that the first expectation on the right-hand side is finite, and we have shown above that the last expectation is also finite. This proves the first claim and completes the proof of  Lemma \ref{lem: entropy}. \hfill{$\Box$}


\subsubsection{Proof of Lemma \ref{lem: Esup}} 
The proof is similar to the proof of Proposition 14 in \cite{azadkia2022linear}, except that it uses Lemma \ref{lem: entropy} above, instead of Proposition 11 of the same paper, thus, we omit the details. \hfill{$\Box$}


\subsubsection{Proof of Lemma \ref{lem: entropy2}} 
We begin with the proof of the first claim.  Lemma \ref{lem: projection} ensures that $\mathcal B_0$ is a compact subset of $\RR^d$ and therefore, we can find $B>0$ such that $\|\beta_0\|< B$ for all $\beta_0\in\mathcal B_0$. This means that the set $\mathcal B_0$ is included in the ball in $\RR^d$ with diameter $B$ in the Euclidean metric, so it follows from Lemma 2.5 in \cite{van2000empirical} that for all $\varepsilon>0$, $\mathcal B_0$ can be covered by $N_\varepsilon:=(4B+\varepsilon)^d\varepsilon^{-d}$ balls of radius $\varepsilon.$ For a fixed $\varepsilon>0$, we denote by $B_j$ and $\beta_j$, with $j\in\{1,\dots,N_\varepsilon\}$ the balls and their centers in such a covering, respectively. We then define the functions $u_j$ and $\ell
_j$ on $\RR^d$ as follows:
\begin{align*}
u_j(x)=\sup_{\beta\in B_j}\beta^\top x \ \ \mbox{and} \ \  \ell_j(x)= \inf_{\beta\in B_j}\beta^\top x.
\end{align*}
Then, for all $\beta\in B_j$, $t< 0$ and $x\in\RR^d$ one has
\begin{align*}
\one_{u_j(x)\leq t}\leq f_{\beta,t}(x)\leq\one_{\ell_j(x)\leq t}.
\end{align*}
This means that the brackets $[\one_{u_j\leq t},\one_{\ell_j\leq t}]$, $j\in\{1,\dots,N_\varepsilon\}$, form a covering of the class $\mathcal F_t$. Similarly, if $t>0$ then  the brackets $[\one_{\ell_j\leq t},\one_{u_j\leq t}]$, $j\in\{1,\dots,N_\varepsilon\}$, form a covering of the class $\mathcal F_t$. We now compute the length of the brackets in the $L_2(\PP^X)$ norm. We have
\begin{align*}
\EE| \one_{u_j(X)\leq t}-\one_{\ell_j(X)\leq t}|^2 
=
\PP(\ell_j(X)\leq t <u_j(X)).
\end{align*}
By the Cauchy-Schwarz inequality,
\begin{align*}
    u_j(x)\leq \beta_j^\top x+\sup_{\beta\in B_j}(\beta-\beta_j)^\top x \leq\beta_j^\top x+\varepsilon\|x\|,
\end{align*}
and similarly, 
$\ell_j(x)\geq \beta_j^\top x-\varepsilon\|x\|$ for all $x$.  Hence,
\begin{align*}
\EE| \one_{u_j(X)\leq t}-\one_{\ell_j(X)\leq t}|^2 
&\leq
\PP(\beta_j^\top X-\varepsilon\|X\|\leq t<\beta_j^\top X+\varepsilon\|X\|)\\
&\leq
\PP(|\beta_j^\top X-t|\leq \varepsilon \|X\|).
\end{align*}
With $C_\varepsilon^p:= \EE\|X\|^p/\varepsilon$, it follows from the Markov inequality that
\begin{align*}
\PP(\|X\|>C_\varepsilon)\leq\varepsilon
\end{align*}
and therefore,
\begin{align*}
\EE| \one_{u_j(X)\leq t}-\one_{\ell_j(X)\leq t}|^2 
&\leq
\PP(|\beta_j^\top X-t|\leq \varepsilon C_\varepsilon) +\varepsilon\\
&=
\mu_0([t-\varepsilon C_\varepsilon,t+\varepsilon C_\varepsilon])+\varepsilon. 
\end{align*}
By assumption there exists $A>0$ such that for all intervals $[a,b]$ in $\RR$, $\mu_0([a,b])\leq A|b-a|$ so we conclude that
\begin{align*}
\EE| \one_{u_j(X)\leq t}-\one_{\ell_j(X)\leq t}|^2 
&\leq
2A\varepsilon C_\varepsilon+\varepsilon\\
&=2A\varepsilon^{1-1/p} \EE^{1/p}\|X\|^p+\varepsilon.
\end{align*}
This means that there exists $A'>0$ (that depends only on $A$ and $\EE^{1/p}\|X\|^p$) such that  
\begin{align*}
\EE^{1/2}| \one_{u_j(X)\leq t}-\one_{\ell_j(X)\leq t}|^2 
\leq \begin{cases}
A'(\varepsilon^{1-1/p})^{1/2}&\mbox{ for all } \varepsilon\in(0,1)\\
A'\varepsilon^{1/2}&\mbox{ for all } \varepsilon>1.
\end{cases}
\end{align*}
This means that the brackets $[\one_{u_j\leq t},\one_{\ell_j\leq t}]$, $j\in\{1,\dots,N_\varepsilon\}$, form a covering of the class $\mathcal F_t$ with length $A'(\varepsilon^{1-1/p})^{1/2}$ if $\varepsilon\in(0,1)$ and $A'\varepsilon^{1/2}$ if $\varepsilon>1$ and therefore,
\begin{align*}
N_{[]}((A'\varepsilon^{1-1/p})^{1/2},\mathcal F_t)
\leq N_\varepsilon
\leq (4B+1)^d\varepsilon^{-d}
\end{align*}
for all $\varepsilon\in(0 ,1)$. Taking $\delta=A'(\varepsilon^{1-1/p})^{1/2}$, this implies that for all $\delta\in(0,A')$,
\begin{align*}
N_{[]}(\delta,\mathcal F_t)
\leq (4B+1)^d(A')^{2dp/(p-1)}\delta^{-2dp/(p-1)}
\end{align*}
whence
\begin{align*}
\log N_{[]}(\delta,\mathcal F_t)
&\leq\log\Big( (4B+1)^d(A')^{2dp/(p-1)}\Big)+\frac{2dp}{p-1}\log\delta^{-1}\\
&\lesssim 1+\log\delta^{-1}.
\end{align*}
Similarly, for $\delta>A'$ we take $\delta=A'\varepsilon^{1/2}$ and obtain that $N_{[]}(\delta,\mathcal F_t)
\leq N_{(\delta/A')^2}$ so that $\log N_{[]}(\delta,\mathcal F_t)\lesssim 1+\log\delta^{-1}$.
The first claim of the lemma follows.

We turn to the proof of the second claim.
It follows from the first claim of the lemma that for all $t\in\RR$ and $\delta>0$,  the bracketing integral of the class $\mathcal F_t$ satisfies
\begin{align*}
{\tilde J}_{[]}(\delta,\mathcal F_t)
&:=\int_0^\delta\sqrt{1+\log N_{[]}(\varepsilon,\mathcal F)}d\varepsilon\\
&\leq
\int_0^\delta\sqrt{1+C+C\log \varepsilon^{-1}}d\varepsilon.
\end{align*}
With arbitrary $q\in(0,2)$, we have $\log\varepsilon^{-1}=(1/q)\log\varepsilon^{-q}\leq\varepsilon^{-q}/q$ and therefore,
\begin{align}\label{eq: majJ}\notag
{\tilde J}_{[]}(\delta,\mathcal F_t)
&\leq \int_0^\delta\sqrt{1+C}+\sqrt{C/q}\varepsilon^{-q/2}d\varepsilon\\
\notag
&=\delta \sqrt{1+C}+\frac {\sqrt {C/q}}{1-q/2}\delta^{1-q/2}\\
&\lesssim\delta+\delta^{1-q/2}.
\end{align}
Now, it follows from  the Cauchy-Schwarz inequality, that for all $t\in\RR$ and $\beta_0\in\mathcal B_0$ one has
\begin{align*}
0\leq f_{\beta_0,t}(x)\leq \one_{B\|x\|>|t|}
\end{align*}
and therefore, the Markov inequality yields that
\begin{align*}
\int f_{\beta_0,t}^2d\PP^X
\leq
\PP(B\|X\|>|t|)
\leq1\wedge\frac{B^p\EE\|X\|^p}{|t|^p}
\end{align*}
where $\EE\|X\|^p$ is finite for $p>2$ by assumption. Now we define $\delta^2=C(1\wedge|t|^{-p})$ for some $C>1$  large enough so that the previous inequality implies that
$$
\int f_{\beta_0,t}^2d\PP^X<\delta^2.
$$
Combining with \eqref{eq: majJ} and Theorem 2.14.17$^\prime$ in \cite{vdvwellner23} proves that 
\begin{align}
\label{eq: 2.14.17'}
\sqrt n\EE\sup_{f\in\mathcal F_t}\Big|\int f(d\PP_n^X-d\PP^X)\Big|
&\lesssim {\tilde J}_{[]}(\delta,\mathcal F_t)\Big(1+\frac{{\tilde J}_{[]}(\delta,\mathcal F_t)}{\delta^2\sqrt n}\Big)
\end{align}
for all $t\in\RR$. Now, consider the case where $|t|>1$. Then by definition, $\delta^2=C|t|^{-p}<C$ and therefore, $\delta< C^{q/4}\delta^{1-q/2}$ which implies that 
$$\delta+\delta^{1-q/2}\lesssim\delta^{1-q/2}\lesssim |t|^{-\frac p2(1-\frac{q}{2})}.$$
Hence, it follows from \eqref{eq: majJ} combined with \eqref{eq: 2.14.17'} that 
\begin{align*}
\sqrt n\EE\sup_{f\in\mathcal F_t}\Big|\int f(d\PP_n^X-d\PP^X)\Big|
&
\lesssim |t|^{-\frac p2(1-\frac q2)}\Big(1+|t|^{\frac p2(1+\frac q2)}n^{-1/2}\Big)\\
&
\lesssim |t|^{-\frac p2(1-\frac q2)}+|t|^{\frac{p q}2}n^{-1/2},
\end{align*}
which completes the proof of the second claim.

To deal with the case where $|t|\leq 1$, we note that by definition, $\delta^2=C$ and $\delta+\delta^{1-q/2}\lesssim 1$. Hence, we can proceed exactly as above just replacing $t$ by $1$ and noting that $n^{-1/2}\lesssim 1$. This completes the proof of Lemma \ref{lem: entropy2}.
\hfill{$\Box$} 

\section*{Acknowledgements}
The research of the third author was supported by the FP2M federation (CNRS FR 2036).

\bibliographystyle{apalike}
\bibliography{bib}
\end{document}